\documentclass[12 pt]{amsart}
\usepackage{amscd,amssymb,amsmath,amsthm}
\usepackage{csquotes}
\usepackage{pdfpages}

\usepackage{tikz}
\usepackage{color}

\usepackage{enumitem}
\usepackage{graphicx}

\newtheorem{Theorem}{Theorem}

\newcommand{\Z}{\mathbb{Z}}

\newcommand{\rarr}{\rightarrow}

\newcommand{\com}{\mathbb{C}}

\newcommand{\ccc}{{\mathrm{ct}}}

\newcommand{\AaA}{\mathsf{A}}
\newcommand{\BbB}{\mathsf{B}}
\newcommand{\M}{{\mathcal{M}}}

\def\E{\mathrm{E}}
\def\n{\mathrm{n}}
\def\L{\mathrm{L}}
\def\V{\mathrm{V}}
\def\H{\mathrm{H}}
\def\g{\mathrm{g}}

\newcommand{\DR}{{\mathsf{DR}}}

\newcommand{\oM}{\overline{\mathcal{M}}}

\newcommand{\cF}{\mathcal{F}}

\newcommand{\AAA}{{\mathsf{A}}}
\newcommand{\BBB}{{\mathsf{B}}}
\newcommand{\CCC}{{\mathsf{H}}}
\newcommand{\B}{\CCC}
\newcommand{\tP}{{\mathcal{P}}}
\newcommand{\cR}{{\mathcal{R}}}
\newcommand{\cS}{{\mathcal{S}}}

\newcommand{\lal}{\left\langle}
\newcommand{\rar}{\right\rangle}

\begin{document}
\baselineskip=16pt

\title{A calculus for the moduli space of curves}

\author{R. Pandharipande}

\date{April 2016}

\maketitle

\setcounter{section}{-1}

\section{Introduction}

The moduli space $\M_g$ of complete nonsingular 
curves of genus $g$ admits a compactification
$$\M_g\subset \overline{\M}_g$$
by stable curves.
Mumford, in  {\em Towards an enumerative geometry of the moduli space of curves} (published in 1983), 
%opened the door to the study of the algebra of tautological classes on $\M_g$.
writes:

\vspace{8pt}
\begin{displayquote}
{{\sf The goal of this paper is to formulate
and to begin an exploration of the enumerative geometry of the set
of all curves of arbitrary genus $g$. By this we mean setting 
up a Chow ring for the moduli space $\M_g$ and its
compactification $\overline{\M}_g$, defining what seem to
be the most important classes in this ring and calculating the
class of some geometrically important loci in $\overline{\M}_g$ in 
terms of these classes. We take as a model for this the
enumerative geometry of the Grassmannians. \cite{Mum}}}
\end{displayquote}

\vspace{8pt}
\noindent Mumford's {\em most important classes} are now termed
{\em tautological classes}. He opened the
door to the study of their algebra --- a fascinating topic
connected to many  areas of modern mathematics.

More than three decades have passed since Mumford's article.
The progress in our understanding of the intersection theory
of the moduli space of curves
has been considerable. Calculations by classical
methods of the 
algebra of tautological classes on $\M_g$ for low $g$ by Faber \cite{Faber}, starting
in the 80's and continuing later in the 90's with Zagier,  have proved
to be fundamental. Witten's conjecture \cite{Wit} in the 90's relating 
the integration of the cotangent line classes on the moduli spaces
$\overline{\M}_{g,n}$ of stable {\em pointed} curves to the KdV hierarchy
was a marvelous  surprise: the study of the
algebra of tautological classes was linked  at a basic
level to the theory of integrable hierarchies. 
The deep role of topology was highlighted in 2007 by the landmark proof 
of Madsen and Weiss \cite{MaW} via homotopy theory
of Mumford's conjecture on the stable cohomology of $\M_{g}$
as $g\rightarrow \infty$.

Starting in the mid 90's,
there was a swift development of Gromov-Witten 
theory.
%{\footnote{See 
%\cite{Beh,BehF,BehM, FulP, KonMan2} for algebraic foundations.}} 
The moduli space $\overline{\M}_{g,n}(X)$ 
of stable maps
intertwines the geometry
of $\overline{\M}_{g,n}$ with the geometry
of the nonsingular target variety $X$. 
Gromov-Witten theory is based upon the virtual fundamental
class \cite{Beh,BehF}
of the moduli of  stable maps,
  $$\Big[\overline{\M}_{g,n}(X)\Big]^{vir} 
\ \in\  A_*(\overline{\M}_{g,n}(X))\, ,$$
a new algebraic cycle{\footnote{All Chow (and cohomology) groups in the paper
will be taken with $\mathbb{Q}$-coefficients.}} whose properties 
constrain the algebra of tautological classes of $\overline{\M}_{g,n}$
in remarkable ways. 

A systematic study of 
the constraints imposed by Gromov-Witten theory
on the algebra of tautological classes was started{\footnote{I have dated
Theorems 1-7 presented in
 the paper  (and the surrounding results) by the years in which the
proofs were found. Published versions appear later and
in mixed order. The dates of publication can be found in the
bibliography.}}
2009 in
\cite{kap} and continued in \cite{PP11,PP13}. Recent progress has culminated in 
a complete proposal by Pixton \cite{Pix} for a calculus of tautological classes on
$\overline{\M}_{g,n}$. 

My goal here is to present 
Pixton's proposal and survey the rapid advances
of the past 6 years. Several open questions are
discussed. An effort has been made to
condense a great deal of mathematics into  as
few pages as possible with the hope that the reader will follow 
through to the end. 

\vspace{12pt}
\noindent{\em \bf Acknowledgments.}
The spirit of my lecture
at the {\em 2015 AMS summer institute in algebraic geometry} in Salt Lake
City has been followed rather closely here. 
I would like to thank the Clay Mathematics Institute for
supporting my visit. 
Some of the material in Sections 1-3 is based on unpublished
notes \cite{joe60} 
of a lecture I gave at {\em A celebration of algebraic geometry}
at Harvard in 2011. Sections 4 and 5 are directly connected to the lecture
in Salt Lake City. A discussion of the recent formula \cite{MOPPZ} for the
Chern characters of the Verlinde bundle on $\overline{\M}_{g,n}$ was presented
in the Salt Lake City lecture, but is omitted here.

Much of what I know about the moduli space of curves has been
learned through collaborations.
Directly relevant to the material presented here is work with 
P. Belorousski, 
C. Faber, G. Farkas, E. Getzler, T. Graber, F. Janda, X. Liu, A. Marian, A. Okounkov, D. Oprea, A. Pixton, 
and D. Zvonkine. Discussions in
 Z\"urich with A. Buryak,
R. Cavalieri,
 E. Clader, D. Petersen, O. Randal-Williams, Y. Ruan, I. Setayesh, and 
Q. Yin have played an important role. 
I have been very fortunate to have had the opportunity to
interact with all of these mathematicians.

 I am supported by the grants
SNF-200020162928 and ERC-2012-AdG-320368-MCSK, and SwissMAP. I am
also supported by  
the Einstein Stiftung in Berlin.

\section{Tautological classes on $\M_g$}

 \subsection{$\kappa$ classes} \label{kcl}
Let $\M_g$ be the moduli space of complete nonsingular
genus $g\geq 2$ curves over $\com$, and let
\begin{equation}\label{xg55}
\pi: \mathcal{C}_g \rightarrow \mathcal{M}_g 
\end{equation}
be the universal curve. We view $\M_g$ and $\mathcal{C}_g$ as
nonsingular, quasi-projective, Deligne-Mumford stacks. 
However, the orbifold
perspective is sufficient for most of our purposes.

The cotangent line $\mathbb{L}$ to the fibers of the morphism \eqref{xg55}
defines a {\em cotangent line class},
$$\psi = c_1(\mathbb{L})\in A^1(\mathcal{C}_g)\ .$$
The $\kappa$ classes are defined by push-forward,
$$\kappa_r = \pi_*( \psi^{r+1}) \in A^{r}(\M_g)\ .$$
The {\em tautological ring} 
$$R^*(\M_g) \subset A^*(\M_g)$$
is the $\mathbb{Q}$-subalgebra generated by all of the
$\kappa$ classes.{\footnote{Since 
$\kappa_{0}= 2g-2 \in \mathbb{Q}$
is a multiple of the fundamental class, we need not take 
$\kappa_0$ as a generator.}}
There is a canonical quotient
$$\mathbb{Q}[\kappa_1,\kappa_2, \kappa_3, \ldots] 
\stackrel{q}{\longrightarrow} R^*(\M_g) \longrightarrow 0\ .$$
The kernel of $q$ is the
ideal of relations among the
$\kappa$ classes.

\subsection{Motivations}
There are two basic motivations for the study of the 
tautological rings $R^*(\M_g)$.
The first is Mumford's conjecture proven
by Madsen and  Weiss \cite{MaW},
$$\lim_{g\rightarrow \infty} H^*(\M_g) \ = \mathbb{Q}[\kappa_1,
\kappa_2, \kappa_3, \ldots ]\, ,$$
determining the {\em stable} cohomology of the moduli 
of curves. While the $\kappa$ classes do not exhaust
$H^*(\M_g)$, there are no other stable classes.
%The perspective of the Faber-Zagier conjecture is to constrain
%the ring of $\kappa$ classes for fixed $g$.

The second motivation comes from a large body of classical
 calculations on $\M_g$ (often related to Brill-Noether
theory). The answers invariably lie in the tautological ring
$R^*(\M_g)$. The study of tautological classes
by Mumford \cite{Mum} was directly inspired  by such
algebro-geometric cycle constructions.

\subsection{Schubert calculus}
The structure of the Chow ring of the Grassmannian $\mathsf{Gr}(r,n)$ of $r$-dimensional
subspaces of  
$\mathbb{C}^n$ is well-known \cite{Ful} and may be viewed as
a model for the study of the tautological classes on $\M_g$. 

The Chern classes of the universal
subbundle
$${\mathsf{S}}\rightarrow \mathsf{Gr}(r,n)\, $$
generate the entire Chow ring,
$$\mathbb{Q} \big[c_1(\mathsf{S}),\ldots,c_r(\mathsf{S})\big]
\stackrel{q}{\longrightarrow} A^*(\mathsf{Gr}(r,n)) \longrightarrow 0\ .$$
The kernel of $q$ is expressed in term of the Segre classes
of the universal subbundle as
$${\text {ker}}(q)= \big(s_{n-r+1}(\mathsf{S}), \ldots,
s_n(\mathsf{S})\big)\, , \ \ \ \ \frac{1}{c(\mathsf{S})}= s(\mathsf{S})\, .$$
The Schubert calculus for the Grassmannian yields 
classical
formulas
for geometric loci in terms of the generators $c_i(\mathsf{S})$.
The subject is fundamentally connected to the representation theory
of the symmetric group.

A basic goal (expressed in the quotation of Mumford in the
Introduction) is to develop a calculus for tautological classes
on the moduli space of curves parallel to the Schubert calculus for the Grassmannian.

\subsection{Cohomology}
We may also define a tautological ring 
$$RH^*(\M_g) \subset 
H^*(\M_g)$$
generated by the $\kappa$ classes in cohomology.
Since there is a natural factoring
$$\mathbb{Q}[\kappa_1,\kappa_2, \kappa_3, \ldots] 
\stackrel{q}{\longrightarrow} R^*(\M_g) 
\stackrel{c}{\longrightarrow} RH^*(\M_g)$$
via the cycle class map $c$, algebraic relations among the
$\kappa$ classes are also cohomological relations.
Whether or not there exist {\em more} cohomological relations is not yet
settled.

\vspace{9pt}
\noindent {\bf Q1.} {\em Is the cycle class map $R^*(\M_g) 
\stackrel{c}{\longrightarrow} RH^*(\M_g)$ an isomorphism?}
\vspace{9pt}

\noindent Calculations (discussed in Section \ref{gorprop} below) show the answer
to question {\bf Q1} is affirmative at least for $g<24$.

\section{Faber-Zagier relations on $\M_g$}

\subsection{Conjecture and proof}
Guided by low genus calculations and deep insight, 
Faber and  Zagier conjectured in 2000 a remarkable set
of relations among the $\kappa$ classes in $R^*(\M_g)$ for all $g$.

The first proof \cite{PP11,PP13} of the Faber-Zagier conjecture (Theorem
\ref{dddd} of Section \ref{3d3d})
 was given in 2010 via
a geometric construction involving the virtual class
of the moduli space of stable quotients \cite{MOP}.{\footnote{The first
relations obtained by virtual localization \cite{GP} on the
moduli space of stable quotients can be found in
\cite[Section 8]{MOP}. See also \cite{Sheng} for a
study of properties of the relations of \cite[Section 8]{MOP}.}}
Since then, a proof \cite{PPZ} in $RH^*(\M_g)$ via Witten's
3-spin class was found in 2013, and a second proof \cite{Janda14,Janda15*}
 in $R^*(\M_g)$ via the virtual class of the moduli space
of stable maps to ${\mathbb{P}}^1$ was found in 2014.

\subsection{Formulas}\label{3d3d}
To write the Faber-Zagier relations, we will require
the following notation.
Let the variable set
$$\mathbf{p} = \{\ p_1,p_3,p_4,p_6,p_7,p_9,p_{10}, \ldots\ \}$$
be indexed by positive integers {\em not} congruent
to $2$ modulo $3$.
Define the series
\begin{multline*}
\Psi(t,\mathbf{p}) =
(1+tp_3+t^2p_6+t^3p_9+\ldots) \sum_{i=0}^\infty \frac{(6i)!}{(3i)!(2i)!} t^i
\\ +(p_1+tp_4+t^2p_7+\ldots) 
\sum_{i=0}^\infty \frac{(6i)!}{(3i)!(2i)!} \frac{6i+1}{6i-1} t^i \ .
\end{multline*}
Since $\Psi$ has constant term 1, we may take the logarithm.
Define the constants $C_r^{\text{\tiny{{\sf FZ}}}}(\sigma)$ by the formula
$$\log(\Psi)= 
\sum_{\sigma}
\sum_{r=0}^\infty C_r^{\text{\tiny{{\sf FZ}}}}(\sigma)\ t^r 
\mathbf{p}^\sigma
\ . $$
The above sum is over all partitions{\footnote{All parts of
a partition are positive (a condition satisfied by the empty
partition).}} $\sigma$ of size 
$|\sigma|$ 
which avoid 
 parts congruent to 2 modulo 3. The empty partition is included
in the sum.
Following standard conventions, to the partition  
$$\sigma=1^{n_1}3^{n_3}4^{n_4} \cdots\, ,$$ we associate
the monomial
$\mathbf{p}^\sigma= p_1^{n_1}p_3^{n_3}p_4^{n_4}\cdots$.
Let 
$$\gamma^{\text{\tiny{{\sf FZ}}}}
= 
\sum_{\sigma}
 \sum_{r=0}^\infty C_r^{\text{\tiny{{\sf FZ}}}}(\sigma)
\ \kappa_r t^r 
\mathbf{p}^\sigma
\ .
$$
For a series $\Theta\in \mathbb{Q}[\kappa][[t,\mathbf{p}]]$ in the variables $
\kappa_i$, $t$, and $p_j$, let
$[\Theta]_{t^r \mathbf{p}^\sigma}$ denote the
 coefficient of $t^r\mathbf{p}^\sigma$
(which is a polynomial in the $\kappa_i$).

\begin{Theorem}[P.-Pixton 2010] \label{dddd} 
{ In $R^d(\M_g)$, the Faber-Zagier relation
$$
\big[ \exp(-\gamma^{\text{\tiny{{\sf FZ}}}}) \big]_{t^d \mathbf{p}^\sigma}  = 0$$
holds when
$d> \frac{g-1+|\sigma|}{3}$ and
$g\equiv d+|\sigma|+1 \mod 2$.}
\end{Theorem}

The dependence upon the genus $g$ in the Faber-Zagier relations of 
Theorem \ref{dddd} occurs in the inequality, the modulo 2 
restriction, and via $\kappa_0=2g-2$. 
For a given genus $g$ and codimension $r$,
Theorem \ref{dddd} provides only {\em finitely} many relations.
While not immediately clear
from the definition,
the $\mathbb{Q}$-linear span of the
Faber-Zagier relations determines an ideal in $\mathbb{Q}[\kappa_1,\kappa_2, \kappa_3, \ldots]$,
see \cite[Section 6]{PP13}.
% --- the matter is discussed in  Section \ref{pppp} and a subset of the Faber-Zagier relations generating the same ideal is described.

\subsection{Gorenstein property}\label{gorprop}
The ring $R^*(\M_g)$ is generated over $\mathbb{Q}$ by the classes
$$\kappa_1,\kappa_2, \ldots, \kappa_{\lfloor \frac{g}{3}\rfloor} \ \in R^*(\M_g)\, ,$$
as conjectured by Faber \cite{Faber}  
and proven by Morita \cite{Morita} in $RH^*(\M_g)$
and Ionel \cite{Ionel} in $R^*(\M_g)$.
By Boldsen's results \cite{Bol},
there are no relations among the $\kappa$ classes of degree less than
or equal to ${\lfloor \frac{g}{3}\rfloor}$.
Looijenga's results \cite{Looij} established 
following nonvanishing and vanishing
conjectures of Faber \cite{Faber}:
\begin{equation}\label{zz12}
R^{g-2}(\M_g) \cong \mathbb{Q}\, , \ \ \ \ R^{>g-2}(\M_g)=0\, .
\end{equation}
The proportionalities in $R^{g-2}(\M_g)$ of monomials
in the $\kappa$ classes are known via Hodge integral
evaluations \cite{Faber, FP1, GetzP}.  The generation, nonvanishing, vanishing,
and proportionality results were all conjectured by Faber in 90s and proven by 2005.

For $g < 24$, the 
Faber-Zagier relations 
yield{\footnote{All the calculations 
discussed in Section \ref{gorprop} concerning
the Faber-Zagier relations are by computer and were undertaken by 
C. Faber in the
period 1991-2011 with continually improving methods. 
The discovery of the failure of $R^*_{\mathsf{FZ}}(\M_{24})$ to be Gorenstein
came in 2009 during a visit to Lisbon.}} 
a Gorenstein
ring,
$$R^*_{\mathsf{FZ}}(\M_g) = \frac{\mathbb{Q}[\kappa_1,\kappa_2,\kappa_3,\ldots]}
{\text{{\sf{FZ}}-relations}}\, ,$$
 with socle in codimension
${g-2}$. By the Gorenstein property,
the pairing
$$R^{d}_{\mathsf{FZ}}(\M_g) \times R^{g-2-d}_{\mathsf{FZ}}(\M_g)\ 
\rightarrow\ R^{g-2}_{\mathsf{FZ}}(\M_g) \cong {\mathbb{Q}}$$
is nondegenerate for $0\leq d \leq g-2$.
The nondegeneracy of the pairing forbids additional relations, so
$$R^*_{\mathsf{FZ}}(\M_g)\cong R^*(\M_g)\, .$$

For $g<24$,
the $\mathbb{Q}$-linear span of the Faber-Zagier relations  
{\em is}  the kernel of
$$\mathbb{Q}[\kappa_1,\kappa_2, \kappa_3, \ldots] 
\stackrel{q}{\longrightarrow} R^*(\M_g) \longrightarrow 0\ $$
 and the cycle class map
$R^*(\M_g) 
\stackrel{c}{\longrightarrow} RH^*(\M_g)$
of question {\bf Q1} {\em is} an isomorphism.

However, the Faber-Zagier relations in genus $24$
do not yield a Gorenstein ring with socle in dimension $22$\,!
There are too few relations in codimension 12,
$$\text{dim}_{\mathbb{Q}}\, R^{12}_{\mathsf{FZ}}(\M_{24}) = 
\text{dim}_{\mathbb{Q}}\, R^{10}_{\mathsf{FZ}}(\M_{24}) + 1\, .$$
Calculations show the Gorenstein property continues
to fail (to  an increasingly greater extent) as $g$ increases above
24. 

\vspace{8pt}
\noindent {\bf Q2.} {\em Do the Faber-Zagier relations
span the ideal of relations among the $\kappa$
classes in $R^*(\M_g)$ for all $g$?}
\vspace{8pt}

While question {\bf Q2} is completely open, a negative answer would
be surprising since many different mathematical approaches 
have failed to find relations outside of the Faber-Zagier
span \cite{Faber, Janda14, PP11, PPZ, RW,  QY}. Moreover, the Gorenstein property for the algebra
of tautological classes has
been {\em proven} to fail for the moduli space   $\M_{2,8}^\ccc$ of
curves of compact type and  the moduli space
$\overline{\M}_{2,20}$ of stable curves
in \cite {Pet,PetT}. So there appears
to be no compelling reason to believe the Gorenstein property 
holds for $\M_{24}$.

\subsection{Hypergeometric series}\label{hyg}
The main actors in the Faber-Zagier relations are the series 
\begin{eqnarray*}
\AaA(t)& = &\sum_{i=0}^\infty \frac{(6i)!}{(3i)!(2i)!} 
%\left(\frac{z}{288}\right)^i\ , \\
t^i\, , \\
\BbB(t)& = & \sum_{i=0}^\infty \frac{(6i)!}{(3i)!(2i)!} \frac{6i+1}{6i-1} 
%\left(\frac{z}{288}\right)^i \ 
t^i \, .
\end{eqnarray*}
%written here in the variable $z=288t$.
The functions $\AaA$ and $\BbB$ are related by the following fundamental
identity observed first by Pixton \cite{Pix}:
\begin{equation}\label{pppp}
 \AaA(-t)\BbB(t)+\AaA(t)\BbB(-t)=-2\, .
\end{equation}

The hypergeometric differential equation satisfied
by $\AaA$, written in the variable $z=288t$, is
$$3 z^2 \frac{d^2}{dz^2} \AaA + (6z-2) \frac{d}{dz} \AaA + \frac{5}{12} \AaA 
= 0 \, .$$
The function $\BbB$ is determined by $\AaA$ and the 
differential equation
$$3 z^2\ \frac{d\AaA}{dz} +\left(\frac{z}{2}-1\right) \AaA = \BbB\, .$$

In the proofs of the Faber-Zagier 
relations, geometric sources for the
series $\AaA$ and $\BbB$ were found. The two approaches 
\cite{Janda14,PP13} in
$R^*(\M_g)$ both find the series in 
the Frobenius geometry associated to $\mathbb{P}^1$.
In the proof of \cite{PPZ} in $RH^*(\M_g)$,
the functions $\AaA$ and $\BbB$ appear in the
$R$-matrix of the Frobenius manifold associated
to $3$-spin curves.

%More recently, occurances
%of $A$ have been noticed \cite{Bur14b}
%in the generating series of descendent integrals over
%the moduli spaces of open Riemann surfaces \cite{PST}. 

\subsection{Descendents}
The series $\AaA$ appears in the 
asymptotic expansion of the Airy function 
related to the Witten-Kontsevich theory
of descendent integration over the moduli space of stable curves.

Let $\overline{\M}_{g,n}$ be 
the moduli space of stable genus $g$ curves with $n$ 
marked points. The cotangent line at the $i^{th}$ marking
defines a line bundle
$$\mathbb{L}_{i} \rarr \overline{\M}_{g,n}\, $$
 with first Chern class
$$\psi_i=c_1(\mathbb{L}_i)\in A^1(\overline{\M}_{g,n})\, .$$ 
The {\em descendent integrals} are defined by
\begin{equation}
\lal\tau_{k_1}\tau_{k_2}\ldots\tau_{k_n}\rar_g=
\int_{\overline{\M}_{g,n}}\psi_1^{k_1}\psi_2^{k_2}\ldots\psi_n^{k_n}\, .
\label{frss}
\end{equation}
The bracket \eqref{frss} vanishes unless the
dimension constraint
$$3g-3+n = \sum_{i=1}^n k_i$$
is satisfied.
The associated  generating series 
in the variables $\{t_i\}_{i\geq 0}$ is
$$
\cF(t_0,t_1,\ldots)=\sum_{\substack{g\geq 0,\,  n\geq 1\\ \vspace{-2pt} \\ 2g-2+n>0}}\frac{1}{n!}\sum_{k_1,\ldots,k_n\ge 0}\lal\tau_{k_1}\tau_{k_2}\ldots\tau_{k_n}\rar_gt_{k_1}t_{k_2}
\ldots t_{k_n}\, .
$$
Witten \cite{Wit} 
conjectured 
%that the 
%partition function $\exp(F)$ is a $\tau$-function of the KdV hierarchy. 
%In particular, 
$$u=\frac{\partial^2 \cF}{\partial t_0^2}$$ is a solution of the KdV hierarchy. 
The first equations of the hierarchy are
\begin{align*}
u_{t_1}&=uu_x+\frac{1}{12}u_{xxx},\\ \notag
u_{t_2}&=\frac{1}{2}u^2u_x+\frac{1}{12}(2u_xu_{xx}+uu_{xxx})
+\frac{1}{240}u_{xxxxx}\, ,\\ \notag
& \vdots
\end{align*}
where we have identified $x$ here with $t_0$. 
Together with the (elementary) string equation,
$$\frac{\partial \cF}{\partial t_0}= \sum_{i=0}^{\infty} t_{i+1} 
\frac{\partial \cF}{\partial t_{i}} + \frac{t_0^2}{2}\,
,$$
Witten's conjecture uniquely determines the series $\cF$.
The conjecture was proven by Kontsevich \cite{Kon92}.
See \cite{KL07,Mir,OP1} for other proofs. 

A basic relationship between the descendent integrals
and the $\AaA$ series is derived via Kontsevich's matrix
integral in \cite{BJP,Kon92}:
\begin{equation}
\label{qqq}
\left.\exp(\cF)\right|_{t_i=-(2i-1)!!\lambda^{-2i-1}}=
\AaA\left(-\frac{\lambda^{-3}}{288}\right)\, .
\end{equation}
For example, the terms of $\exp(\cF)$ which contribute to the
$\lambda^{-3}$ coefficient are:
\begin{eqnarray*}
\lal \tau_0^3 \rar_0 \frac{\left(-\lambda^{-1}\right)^{3}}{3!} + \lal \tau_1 \rar_1
\left(-\lambda^{-3}\right) &= &\left(\frac{1}{6}+\frac{1}{24}\right) 
\left(-\lambda^{-3}\right) \\
& = & \frac{1440}{24} \left(-\frac{\lambda^{-3}}{288}\right) \\
& = & \frac{6!}{3!2!} \left(-\frac{\lambda^{-3}}{288}\right)\, .
\end{eqnarray*}
The result is the $\lambda^{-3}$ term of 
$\AaA\left(-\frac{\lambda^{-3}}{288}\right)$.

\vspace{8pt}
\noindent {\bf Q3.} {\em Is the descendent evaluation \eqref{qqq}
related to the occurrence of the $\AaA$ series in the relations among
tautological classes?}
\vspace{8pt}

Question {\bf Q3} was first raised in \cite{BJP} where 
the many appearances
of the $\AaA$ and $\BbB$ series related to the moduli spaces of curves
are surveyed.

\subsection{Markings}
The algebra of tautological classes may also be studied on the moduli spaces
$${\M}_{g,n}\, , \ \ {\M}_{g,n}^{\mathrm{rt}}\, , \ \ \mathcal{C}_g^n\, $$
of {\em pointed} curves lying over $\M_g$. 
Here, $\M_{g,n}$ is the moduli space of nonsingular curves with distinct
markings,
${\M}_{g,n}^{\mathrm{rt}}$ is the moduli of curves with rational
tails, and $\mathcal{C}_g^n$ is the $n^{th}$ fiber product of the universal curve
$$\pi: \mathcal{C}_g \rightarrow \M_g\, .$$
In addition to the $\kappa$ classes, there are cotangent line classes
$\psi_i$ at the markings and  diagonal classes.
The full algebra of tautological classes 
will be considered in Section \ref{ll99}.

%\subsection{Boundary}
%As a corollary of our proof of Theorem \ref{dddd} via the moduli space of stable
%quotients, we obtain the following stronger boundary result.
%If $g-1+|\sigma|< 3r$ and
%$g\equiv r+|\sigma|+1 \mod 2$, then
%\begin{equation}
%\big[ \exp(-\gamma^{\text{\tiny{{\sf FZ}}}}) \big]_{t^r \mathbf{p}^\sigma}  \in% R^*(\partial\overline{\mathcal{M}}_g)\ 
%\end{equation}
%Not only is the Faber-Zagier relation 0 on $R^*(\mathcal{M}_g)$, but the
%relation is equal to a tautological class on the boundary of
%the moduli space $\overline{\mathcal{M}}_g$. A precise conjecture for 
%the boundary terms has been proposed in \cite{Pix}.

\section{$\kappa$ classes on $\M_{g,n}^\ccc$}
\subsection{Compact type}
A connected nodal curve $C$ is of {\em compact type} if the dual graph of $C$
is a tree (or, equivalently, if the Jacobian of $C$ is compact). 
The moduli space $\M_{g,n}^{\mathrm{ct}}$ of curves of compact type
lies in between
the moduli spaces of nonsingular and stable curves
$$\M_{g,n} \, \subset\, \M_{g,n}^{\mathrm{ct}} \, \subset\, \overline{\M}_{g,n}\, .$$

The algebra of $\kappa$ classes on $\M_{g,n}^{\mathrm{ct}}$
is remarkably well-behaved. The structure is simpler and
better understood than for $\M_g$.
Of course,
the full algebra of tautological classes on $\M_{g,n}^{\mathrm{ct}}$
contains more than
the $\kappa$ classes (and will be discussed in Section \ref{ll99}).

\subsection{$\kappa$ classes}
The definition of the $\kappa$ classes on $\M_g$ is easily extended to
the moduli of stable curves. 
%The $\kappa$ classes in the Chow ring 
%$A^*(\overline{\M}_{g,n})$
%are defined by the following geometry.
Let
$$\pi: \overline{\M}_{g,n+1} \rarr \overline{\M}_{g,n}$$
be the universal curve viewed as the ($n+1$)-pointed space, and let
$$\psi_{n+1} = c_1(\mathbb{L}_{n+1})\ \in A^1(\overline{\M}_{g,n+1})$$
be the Chern class of the cotangent line at the last marking.
The $\kappa$ classes 
are
$$\kappa_r = \pi_*(\psi^{r+1}_{n+1})\  \in A^r(\overline{\M}_{g,n}), 
\ \ \ i \geq 0\ .$$
The simplest is $\kappa_0$ which equals $2g-2+n$ times the unit
in $A^0(\overline{\M}_{g,n})$.
%The convention 
%$$\kappa_{-1}= \epsilon_*(\psi_{n+1}^0)=0$$
%is often convenient. 

The restriction of $\kappa_r\in A^r(\overline{\M}_g)$ to
$\M_g$ via the inclusion
$$\M_g \, \subset \, \overline{\M}_g$$
agrees with the definition of $\kappa_r\in A^r(\M_g)$
of Section \ref{kcl}.
The $\kappa$ classes on $\M_{g,n}^\ccc$
are defined via restriction from $\overline{\M}_{g,n}$.

Define the $\kappa$ ring for curves of compact type,
$$
\kappa^*(\M_{g,n}^{\mathrm{ct}})\, \subset \, A^*(\M_{g,n}^\ccc)\, , 
$$
to be the $\mathbb{Q}$-subalgebra 
generated by the $\kappa$ classes.
Of course,  the $\kappa$ rings are graded by degree.

\subsection{Basic results}

%\vspace{15pt}
%\noindent {\bf D. Universality for $M_{g,n}^c$}
%\vspace{15pt}

The ring $\kappa^*(\M^\ccc_{g,n})$ is generated over $\mathbb{Q}$ by the classes
$$\kappa_1,\kappa_2, \ldots, \kappa_{g-1+\lfloor\frac{n}{2}
\rfloor} \ \in \kappa^*(\M_{g,n}^\ccc)\, .$$
If $n>0$, there are no relations of degree less than or equal to
${g-1+\lfloor\frac{n}{2}
\rfloor}$. In $\kappa^*(\M_g^\ccc)$, there are no relations of degree
less than $g-1$ (whether degree $g-1$ relations can occur is
not known).
The proofs of the above generation and freeness results
can be found in \cite{kap}.

The nonvanishing
and vanishing results for $\kappa^*(\M_{g,n}^\ccc)$, 
$$\kappa^{2g-3+n}(\M^\ccc_{g,n}) \cong \mathbb{Q}, \ \ \ \ \
\kappa^{>2g-3+n}(\M^\ccc_{g,n})=0\, , $$
parallel to \eqref{zz12}
for the $\kappa$ classes on $\M_g$,
were proven in \cite{FPrel,GV}. The proportionalities
in $\kappa^{2g-3+n}(\M^\ccc_{g,n})$ of the monomials in the $\kappa$
classes are determined by Hodge integral evaluations \cite{FP1, FP2, GetzP}.

\subsection{Universality}

A new and surprising feature about the $\kappa$ rings in the compact
type case is the following universality result proven in \cite{kap}.
There appears to be no parallel for $\M_g$.

\begin{Theorem}[P. 2009] \label{dxxd}
Let $g>0$ and $n>0$, then the assignment $\kappa_i \mapsto \kappa_i$
extends to a ring isomorphism
$$ \iota: \kappa^*(\M_{g-1,n+2}^{\ccc}) \cong  \kappa^*(\M_{g,n}^\ccc) \ . $$
\end{Theorem}

In other words, the relations among the $\kappa$ classes
in the above cases 
are genus independent! 
By composing the isomorphisms $\iota$ of Theorem 2, we obtain isomorphisms
$$\iota: \kappa^*(\M_{0,2g+n}^\ccc) \cong \kappa^*(\M_{g,n}^\ccc)$$
so long as $n>0$. Hence, universality reduces all questions about
the $\kappa$ rings to genus 0. 
Calculations of the relations,
bases, and Betti numbers of the ring
$\kappa^*(\M_{g,n>0}^\ccc)$ are 
obtained in \cite{kap} using the genus 0 reduction.

Let $P(d)$ be the set of partitions of $d$, and let
$$P(d,k)\subset P(d)$$ be the set of partitions of $d$ into
at most $k$ parts. 
Let $|P(d,k)|$ be the
cardinality. To a partition with positive parts
$\mathbf{p}= (p_1,\ldots,p_\ell)$ in $P(d,k)$,
we associate a $\kappa$ monomial by
$$\kappa_{\mathbf{p}} = \kappa_{p_1} \cdots \kappa_{p_\ell} \in 
\kappa^d(\M_{g,n}^\ccc) \ .$$

\begin{Theorem}[P. 2009] \label{htt5} For $n>0$,
a $\mathbb{Q}$-basis of $\kappa^d(\M_{g,n}^\ccc)$ is given by
$$\{ \kappa_{\mathbf{p}} \ | \ \mathbf{p} \in P(d,2g-2+n-d)\ \} \  .$$
\end{Theorem}

The Betti number calculation,
$$\text{dim}_{\mathbb{Q}} \ \kappa^d(\M_{g,n}^\ccc) \ = \ |P(d,2g-2+n-d)|\ ,$$
is implied by Theorem \ref{htt5}.

\vspace{8pt}
\noindent {\bf Q4.} {\em Is there a representation theoretic
formula for the multiplication of $\kappa^*(\M_{g,n}^\ccc)$
in the basis of Theorem \ref{htt5}?}
\vspace{8pt}

Question {\bf Q4} appears directly parallel to
questions in the Schubert calculus. Recent progress on {\bf Q4}
has been made by Setayesh \cite{Set} who has found a formula
involving the combinatorics of partitions.

Theorems \ref{dxxd} and \ref{htt5}  for the $\kappa$ rings in the compact type
case require at least 1 marked point. In the unmarked case,
half of the universality still holds. The
assignment $\kappa_i\mapsto \kappa_i$ extends to a surjection
$$\iota_g: \kappa^*(\M_{0,2g}^\ccc) \rightarrow \kappa^*(\M_{g}^\ccc)\ . $$
However, a nontrivial kernel is possible. The first kernel
occurs in genus $g=5$.

\vspace{+8pt}
\noindent{\bf Q5.} {\em What is the kernel of
$\kappa^*(\M_{0,2g}^\ccc) \rightarrow \kappa^*(\M_{g}^\ccc)$?}
\vspace{+8pt}

In genus 5, the kernel of $\iota_5$ is related to Getzler's
relation in $\overline{\M}_{1,4}$, see \cite{kap} for a discussion.
A complete answer to {\bf Q5}  will likely involve sequences
of special relations in the
tautological ring.{\footnote{ {\bf Q5} is perhaps
the narrowest question discussed in the paper, but I am very
curious to know the answer.}}

The proofs of Theorems \ref{dxxd} and \ref{htt5} were obtained in \cite{kap} by studying the
virtual class of the moduli space of stable quotients \cite{MOP}. The results were
precursors to the proof \cite{PP11,PP13} 
of the Faber-Zagier relations via the moduli of stable quotients.

\subsection{Relations}
A natural question to ask is whether $\kappa$ relations
in the compact type case can be put in a  form parallel to Theorem 1.
An answer has been found by Pixton \cite[Section 3.4]{PixPhd}.

We define a
set of relations as follows.
Let 
$$\mathbf{p} = \{\ p_1,p_2,p_3,\ldots\ \}$$
be a variable set indexed by all positive integers.
Let
\begin{multline*}
\ \ \Psi(t,\mathbf{p}) =
(1+tp_2+t^2p_4+t^3p_6+\ldots) \sum_{i=0}^\infty (2i-1)!!\ t^i
\\ +(p_1+tp_3+t^2p_5+\ldots)\, , \  \ \ \ \ \ \ \ \ \ \ \ \ \ \ \ \ \ \ 
\end{multline*}
where $(2i-1)!!= \frac{(2i)!}{2^ii!}$ as usual.
Define the constants $C_r^{\mathsf{P}}(\sigma)$ by the formula
$$\log(\Psi)= 
\sum_{\sigma}
\sum_{r=0}^\infty C_r^{\mathsf{P}}(\sigma)\ t^r 
\mathbf{p}^\sigma
\ . $$
Here, $\sigma$ denotes any partition with positive parts.
Let 
$$\gamma^{\mathsf{P}}= 
\sum_{\sigma}
 \sum_{r=0}^\infty C_r^{\mathsf{P}}(\sigma)
\ \kappa_r t^r 
\mathbf{p}^\sigma
\ .
$$

\begin{Theorem}[Pixton 2013] \label{ctp}
In $\kappa^d(\M_{g,n}^\ccc)$, the relation
$$
\big[ \exp(-\gamma^{\mathsf{P}}) \big]_{t^d \mathbf{p}^\sigma}  = 0$$
holds when
$d>\frac{2g-2+n+|\sigma|}{2}$.
\end{Theorem}

Theorems \ref{dxxd} and \ref{htt5} determine the complete set
of relations among the $\kappa$ classes in the $n>0$ case.
Using Theorems \ref{dxxd} and \ref{htt5},
Pixton \cite{PixPhd} has proven that the
$\mathbb{Q}$-linear span of the 
relations of Theorem \ref{ctp} {\em is} the complete set
in the $n>0$ case.

\section{Pixton's relations on $\overline{\M}_{g,n}$} 
\label{Ssec:relations}

\subsection{Overview}
Tautological classes on 
the moduli space $\overline{\M}_{g,n}$ of stable curves are obtained
from $\kappa$ classes, $\psi$ classes, {\em and}
the classes of boundary strata (indexed by stable graphs).{\footnote{A
study of the algebra of $\kappa$ classes can also be pursued
on $\oM_{g,n}$, see \cite{Set1,Set2}.}}
Decorated stable graphs provide a language for describing
all tautological classes. A parallel role is played by the
 language of partitions
in  the Schubert calculus of the
Grassmannian.
After a brief discussion of stable graphs in Section \ref{bs}, the algebra of tautological classes
$$R^*(\overline{\M}_{g,n}) \subset A^*(\overline{\M}_{g,n})$$
is defined in Section \ref{straa}. 

The main recent advance is the set of relations in 
$R^*(\overline{\M}_{g,n})$
found by Pixton.
 Pixton's relations \cite{Pix} were
 conjectured in 2012 
and first proven \cite{PPZ} to hold in 
$RH^*(\overline{\M}_{g,n})$ in 2013
using Witten's 3-spin class{\footnote{Tautological relations obtained from
Witten's $r$-spin class for higher $r$ are studied in \cite{PPZ16}. By
Janda's result \cite{Janda15}, the relations of \cite{PPZ16} are contained in Pixton's set.}} 
 and the
Givental-Teleman classification of Cohomological Field Theories \cite{Tel}. 
Shortly afterwards (also in 2013), Janda \cite{Janda13}   
found a proof in $R^*(\overline{\M}_{g,n})$ 
using a mix of virtual localization \cite{GP} and
$R$-matrix techniques for the equivariant stable
quotients theory of $\mathbb{P}^1$. His argument combined
elements of the Chow results of \cite{PP11,PP13} and
the CohFT methods of \cite{PPZ}. A second proof in Chow
via the equivariant Gromov-Witten theory of $\mathbb{P}^1$
was found by Janda \cite{Janda15*} in 2015.

Pixton has further conjectured that  {\em all}  
relations among tautological classes
are obtained from his set. The claim
is open and is perhaps
the most important question in the subject.
Pixton's relations for $R^*(\overline{\M}_{g,n})$ specialize to the
Faber-Zagier relations when restricted to the moduli space
${\M}_g$ of nonsingular curves.

%The 
%Chow results of \cite{Jan13} are via a study of  the equivariant Gromov-Witten
%theory of $\mathbb{CP}^1$. The application of
%the virtual localization formula \cite{GrP} to the Gromov-Witten theory
%of $\mathbb{CP}^1$ bypasses the need for Teleman's cohomological
%results, and the proof of \cite{Jan13} holds in Chow.

%Both Witten's 3-spin theory and the equivariant Gromov-Witten
%theory of $\mathbb{CP}^1$ are {\em Chow Field Theories}: the CohFT
%axioms are satisfied in the Chow rings $A^*(\oM_{g,n},\mathbb{Q})$.

%\begin{question}
%  Does the Givental-Teleman classification hold for semisimple Chow
%  Field theories?
%\end{question}

%The proof of \cite{PPZ13} would imply the Chow vanishing of
%Theorem \ref{Thm:relations} if the answer to Question 4 is yes.
%The classification of CohFTs by Teleman uses the stable
%cohomology of the moduli spaces of curves. Are such stability
%results valid in the Chow ring?

\subsection{Boundary strata} \label{bs}
The boundary
strata of the moduli space $\overline{\M}_{g,n}$  
parameterizing complex structures on curves of
 {\em fixed topological type} correspond
to {\em stable graphs}. The idea here is simple,
but the notation requires some care.

A {\em stable graph} $\Gamma$ 
consists of the data
$$\Gamma=(\V, \H,\L, \ \mathrm{g}:\V \rarr \Z_{\geq 0},
\ v:\H\rarr \V, 
\ \iota : \H\rarr \H)$$
which satisfies the following properties:
\begin{enumerate}
\item[(i)] $\V$ is a vertex set with a genus function 
$\g:\V\to \Z_{\geq 0}$,
\item[(ii)] $\H$ is a half-edge set equipped with a 
vertex assignment $$v:\H \to \V$$ and an 
involution $\iota:\H\to \H$,
\item[(iii)] $\E$, the edge set, is defined by the
2-cycles of $\iota$ in $\H$ (self-edges at vertices
are permitted),
\item[(iv)] $\L$, the set of legs, is defined by the fixed points of $\iota$ and is endowed with a bijective correspondence with the set of markings $$\L \leftrightarrow \{1,\ldots, n\}\, ,$$
\item[(v)] the pair $(\V,\E)$ defines a {\em connected} graph,
\item[(vi)] for each vertex $v$, the stability condition holds:
$$2\g(v)-2+ \n(v) >0,$$
where $\n(v)$ is the valence of $\Gamma$ at $v$ including 
both edges and legs.
\end{enumerate}
An automorphism of $\Gamma$ consists of automorphisms
of the sets $\V$ and $\H$ which leave invariant the
structures $\mathrm{g}$, $\iota$, and $v$ (and hence respect $\E$ and $\L$).
Let $\text{Aut}(\Gamma)$ denote the automorphism group of $\Gamma$.

The genus of a stable graph $\Gamma$ is defined by
$$\g(\Gamma)= \sum_{v\in V} \g(v) + h^1(\Gamma).$$
A boundary stratum of the moduli space $\oM_{g,n}$ 
 naturally determines
a stable graph of genus $g$ with $n$ legs by considering the dual graph of a generic pointed curve parameterized by the stratum.

To each stable graph $\Gamma$, we associate the moduli space
\begin{equation*}
\oM_\Gamma =\prod_{v\in \V} \oM_{\g(v),\n(v)}.
\end{equation*}
 Let $\pi_v$ denote the projection from $\oM_\Gamma$ to 
$\oM_{\g(v),\n(v)}$ associated to the vertex~$v$.  There is a
canonical
morphism 
\begin{equation}\label{dwwd}
\xi_{\Gamma}: \oM_{\Gamma} \rarr \oM_{g,n}
\end{equation}
 with image{\footnote{
The degree of $\xi_\Gamma$ is $|\text{Aut}(\Gamma)|$.}}
equal to the closure of the boundary stratum
associated to the graph $\Gamma$.  To construct $\xi_\Gamma$, 
a family of stable pointed curves over $\oM_\Gamma$ is required.  Such a family
is easily defined 
by attaching the pull-backs of the universal families over each of the 
$\oM_{\g(v),\n(v)}$  along the sections corresponding to half-edges.
Let 
$$[{\Gamma}] \in A^*(
{\overline{\mathcal{M}}_{g,n}})$$
denote the push-forward under $\xi_\Gamma$ of the fundamental
class of $\oM_{\Gamma}$.

Two examples of boundary strata in $\overline{\M}_{3,3}$ and their associated stable graphs are given in the following diagram. 

{\centering{\ \ \ \ \hspace{50pt} \includegraphics[width=10cm]{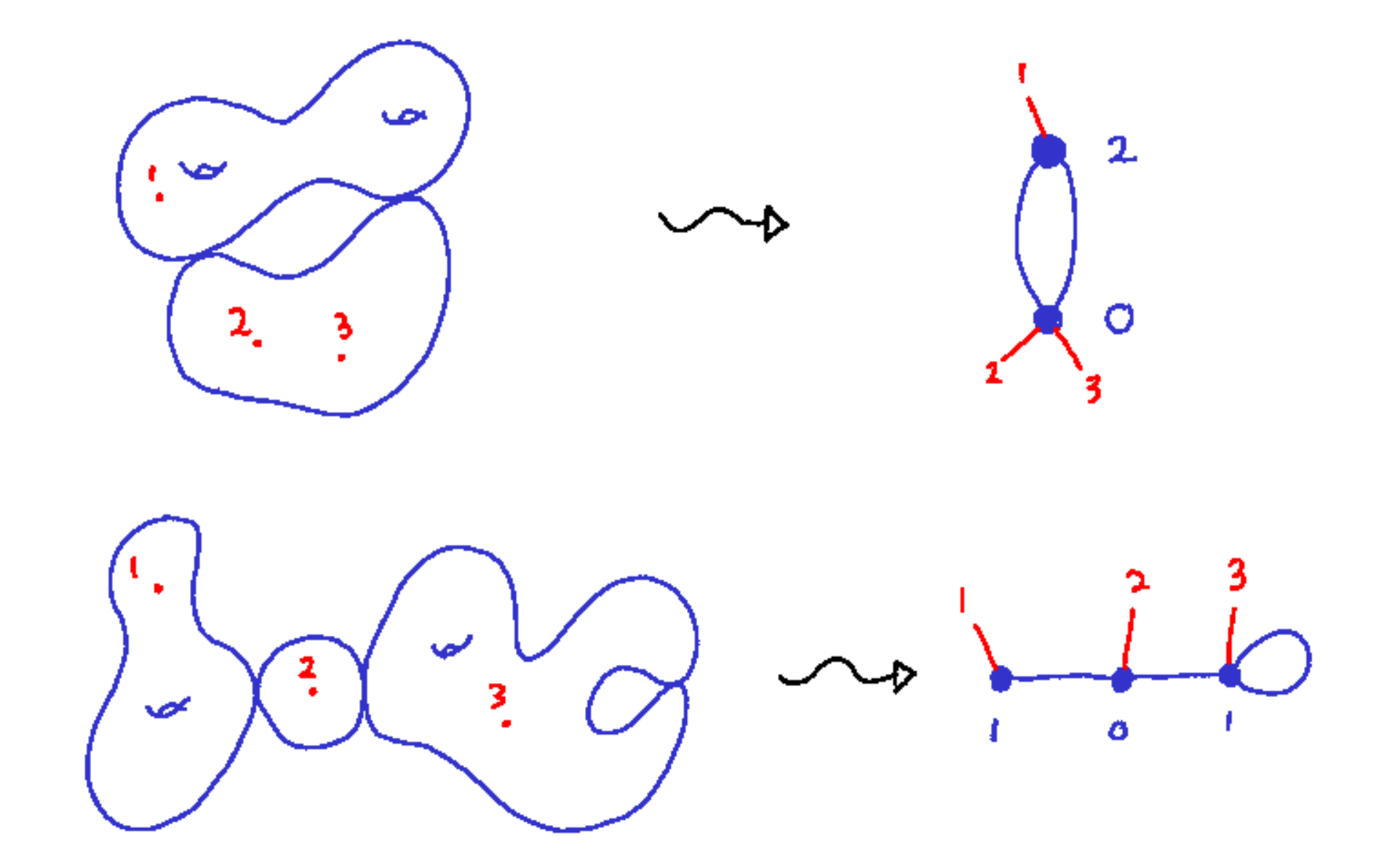}}}

\noindent Let $\Phi$ and $\widehat{\Phi}$ be the stable graphs
in the first and second cases in the diagram.
The moduli space $\oM_{\Phi}$ in the first case is
$$\oM_{0,4} \times\oM_{2,3}$$
with $|\text{Aut}(\Phi)|=2$. The nontrivial automorphism
arises from switching the edges.
The moduli space $\oM_{\widehat{\Phi}}$ in the second case is 
$$\oM_{1,2} \times \oM_{0,3} \times \oM_{1,4}$$
with $|\text{Aut}(\widehat{\Phi})|=2$. The nontrivial automorphism
arises from switching the half-edges on the
self-edge.

\subsection{Strata algebra} \label{straa}
Let $\cS_{g,n}^*$ be the $\mathbb{Q}$-algebra of
$\kappa$ and $\psi$ classes 
supported on the strata of $\overline{\mathcal{M}}_{g,n}$.
A $\mathbb{Q}$-basis of $\cS_{g,n}^*$ is given by isomorphism classes
of pairs $[\Gamma, \gamma]$ where $\Gamma$ is a 
stable graph corresponding to a stratum of the moduli space, 
$$\overline{\mathcal{M}}_\Gamma \rightarrow \overline{\mathcal{M}}_{g,n}\, ,$$
and $\gamma$ is a product of $\kappa$ and $\psi$
classes on $\overline{\mathcal{M}}_\Gamma$. The $\kappa$ classes
are associated to the vertices, and the $\psi$ classes are
associated to the half-edges. The only condition imposed is
that the degrees of the $\kappa$ and $\psi$ classes 
associated to a vertex $v\in \V(\Gamma)$ together do {\em not} exceed
the dimension 
$3\g(v)-3+\n(v)$ of the moduli space at $v$.

For the graph $\Phi$ associated to a stratum of $\overline{\M}_{3,3}$
in the diagram, let $v_0$ and $v_2$ denote the vertices of genus $0$ and
$2$ respectively. Let the left edge consist of the  two half edges
$h_0-h_2$ where $h_0$ is incident to $v_0$ and $h_2$ is incident to $v_2$.
Then,
\begin{equation}\label{f334}
\left[\Phi, \kappa_1[v_0]\kappa_2[v_2]\psi_{h_2}^2\psi_1\right]
\end{equation}
is an example of such a pair. The codimension of the pair
\eqref{f334} is 
$$8=2+6\, ,$$
 2 for the
nodes of $\Phi$ and 6 for the $\kappa$ and $\psi$ classes.

The strata algebra $\cS_{g,n}^*$ is graded by
codimension
$$\cS_{g,n}^* = \bigoplus_{d=0}^{3g-3+n} \cS^d_{g,n}\ $$
and carries a product for which
the natural push-forward map
\begin{equation}\label{v123}
\cS_{g,n}^* \rightarrow A^*(\overline{\mathcal{M}}_{g,n})
\end{equation}
is a ring homomorphism,
 see \cite[Section 0.3]{PPZ} 
for a detailed discussion.

The image of \eqref{v123} is, by definition, the {\em tautological ring}{\footnote{Our definition here follows the Appendix of \cite{GP2}.
See \cite[Section 1]{FP3} for a more intrinsic approach.}}
$$R^*(\overline{\mathcal{M}}_{g,n})\subset A^*(\overline{\mathcal{M}}_{g,n})\ .$$
Hence, we have a quotient
$$ \cS_{g,n}^*
\stackrel{q}{\longrightarrow} R^*(\overline{\mathcal{M}}_{g,n}) \longrightarrow 0\ .$$
The ideal of {\em tautological relations} is the kernel of $q$.

In the strata algebra, the basis elements $[\Gamma,\gamma]$
are treated formally.  In the case $(g,n)=(0,4)$, we have 
$$\text{dim}_{\mathbb{Q}}\, \cS^0_{0,4}= 1\, , \ \ \
\text{dim}_{\mathbb{Q}}\, \cS^1_{0,4}= 8\, .$$
Let $\Gamma_{\bullet}$ be the unique graph of genus 0 with 4 markings 
and a single vertex $v$. Then $$[\Gamma_{\bullet},1]\in \cS^0_{0,4}\, $$
is a basis.
In codimension 1,
 the
5 possibilities for $\gamma$ on $\Gamma_{\bullet}$ yield the pairs
$$[\Gamma_\bullet,\kappa_1[v]]\, , \ 
[\Gamma_\bullet,\psi_1]\, , \ 
[\Gamma_\bullet,\psi_2]\, , \ 
[\Gamma_\bullet,\psi_3]\, ,\  
[\Gamma_\bullet,\psi_4]\, \in \cS^1_{0,4} .$$
In addition, there are 3 pairs
$$[\Gamma_{1,2|3,4},1]\, , \  [\Gamma_{1,3|2,4},1]\, ,  \ [\Gamma_{1,4|2,3},1]\, 
\in \cS^1_{0,4}$$ 
where the underlying graphs have 2 vertices (and correspond
to the usual boundary strata).
The kernel 
$$\cS^1_{0,4} \stackrel{q}{\longrightarrow} R^1(\overline{\mathcal{M}}_{0,4})$$
is 7 dimensional and contains
the basic 
linear equivalence of the three boundary divisors of 
$$\oM_{0,4}\cong {\mathbb{P}}^1\, ,$$
see \cite{keel,KonMan} for a study of $R^*(\oM_{0,n})$.

The first
 geometrically interesting relation{\footnote{A proof of Getzler's
relation in Chow was given later in \cite{PPP}.}} was found in genus 1 by Getzler 
\cite{Get} in 1996.
Soon after, several low genus relations were determined. Below is 
a tautological relation{\footnote{The strata classes
in the genus 2 relation of the diagram have been represented
by their {\em topological type} instead of their associated dual graph.
The genera of the components are underlined. The red marked points are
unlabeled. Each picture represents the sum of the 6 possible
labelings of the markings. The diagram (taken from \cite{BP}) 
was typeset by P. Belorousski.}}
  in codimension 2 on $\oM_{2,3}$ found in \cite{BP} in
1998.

\vspace{5pt}

{{ \hspace{-18pt} \includegraphics[width=13.2cm]{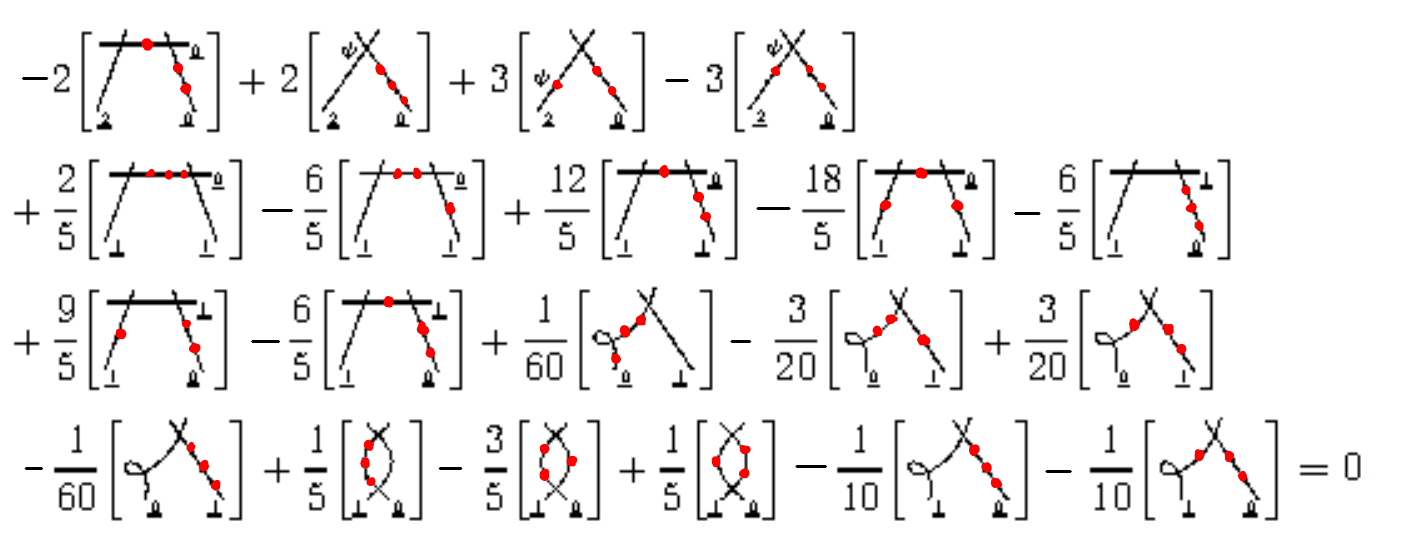}}}

\vspace{5pt}

\noindent Pixton's set puts order to the chaos of strata classes and coefficients
which appear in the above relation (and in all the other relations{\footnote{For example, see
\cite{Getz22,Kim1,Kim2} for further relations in genus 2 and 3.}} 
found
in the period after Getzler's discovery).

\subsection{Vertex, leg, and edge factors} 
Pixton's relations are determined by 
a set 
$$\tP=\{\cR_{g,A}^d\}$$ of
elements $\cR^d_{g,A}\in \cS^d_{g,n}$ associated to the data
\begin{enumerate}
\vspace{7pt}
\item[$\bullet$] $g,n\in \mathbb{Z}_{\geq 0}$ in the stable range $2g-2+n>0$,
\vspace{7pt}
\item[$\bullet$] $A=(a_1,\ldots, a_n), \ \, a_i \in\{0,1\}$,
\vspace{7pt}
\item[$\bullet$] $d\in \mathbb{Z}_{\geq 0}$ satisfying
$d > \frac{g-1+\sum_{i=1}^n a_i}{3}$.
\end{enumerate}
\vspace{7pt}
The elements $\cR^d_{g,A}$ are expressed as sums over
stable graphs of genus $g$ with $n$ legs. 
Before writing the formula for $\mathcal{R}^d_{g,A}$, a few
 definitions are required.

\vspace{9pt}
\noindent {\bf Definitions of  $\CCC_0$ and $\CCC_1$.} The
hypergeometric series $\AAA$ and $\BBB$ of Section \ref{hyg} enter
Pixton's relations in following form:

{\small
\begin{align*}
\CCC_0(T)&=\AAA(-T) = \sum_{i=0}^\infty \frac{(6i)!}{(2i)!(3i)!}(-T)^i
= 1-60T+27720T^2 -\cdots,\\ 
\CCC_1(T)&=-\BBB(-T) = -\sum_{i= 0}^\infty  \frac{(6i)!}{(2i)!(3i)!}\frac{6i+1}{6i-1}
(-T)^i
=1 + 84T - 32760T^2 + \cdots.
\end{align*}}
%\noindent These series control the original Faber-Zagier relations
% and continue
%to play a central role in the set $\tP$.

\vspace{4pt}
\noindent{\bf Definition of $\kappa(f)$.}
Let $f(T)$ be a power series with vanishing constant and linear terms,
$$f(T)\in T^2\mathbb{Q}[[T]]\ .$$
For each $\oM_{g,n}$,  we define
\begin{equation}\label{g33g}
\kappa(f) = \sum_{m \geq 0} \frac1{m!}\ { p_{m*}} \Big(f(\psi_{n+1}) \cdots f(\psi_{n+m})\Big)
 \ \in A^*(\oM_{g,n})\, ,
\end{equation}
where $p_m$ is the forgetful map
$$p_m: \oM_{g,n+m} \to \oM_{g,n}\, .$$ 
By the vanishing in degrees 0 and 1  of $f$, the sum \eqref{g33g} is finite.

\vspace{12pt}
\noindent{\bf Definitions of $\mathsf{G}_{g,n}$ and $\zeta_v$.}
Let $\mathsf{G}_{g,n}$ be the (finite) set of stable graphs of
genus $g$ with $n$ legs (up to isomorphism).
Let $\Gamma \in \mathsf{G}_{g,n}$. For each vertex $v\in \V$,
we introduce an auxiliary variable $\zeta_v$ and impose the
conditions
$$\zeta_v \zeta_{v'}= \zeta_{v'} \zeta_v\, , \ \ \ \zeta_v^2= 1\, .$$
 The variables $\zeta_v$ will be responsible for keeping track of
a local parity condition at each vertex.

\vspace{12pt}
The formula for $\mathcal{R}_{g,A}^d$ is a sum over $\mathsf{G}_{g,n}$.
The summand corresponding to $\Gamma \in \mathsf{G}_{g,n}$ is a 
product of 
vertex, leg, and edge factors:
\begin{enumerate}
%\vspace{5pt}
\item[$\bullet$]
For $v\in \V$, let
$\kappa_v = \kappa\big(T-T \B_0(\zeta_vT)\big)$.
\vspace{8pt}
\item[$\bullet$]
For $l \in \L$,
let
 $\B_l =\zeta_{v(l)}^{a_l} \B_{a_l} \! \left(\zeta_{v(l)} \psi_{l}\right)$,  
where $v(l)\in V$ is the vertex to which the leg is assigned.
\vspace{8pt}
\item[$\bullet$]
For $e\in \E$, let
\vspace{5pt}
\begin{align*}
\Delta_e &= \frac{ \zeta' + \zeta'' - 
\B_0(\zeta' \psi') \zeta''\B_1(\zeta'' \psi'')
-\zeta'\B_1(\zeta' \psi') \B_0(\zeta'' \psi'')}
{\psi'+\psi''} \\ 
& \vspace{5pt} =(60 \zeta' \zeta''-84) +
\left[32760(\zeta'\psi' + \zeta'' \psi'') - 27720 (\zeta'\psi''+\zeta''\psi')\right] + \cdots,
\end{align*}

\vspace{5pt}
\noindent where $\zeta',\zeta''$ are the $\zeta$-variables assigned to the vertices adjacent to the edge $e$ and $\psi', \psi''$ are the $\psi$-classes corresponding to the half-edges.
\end{enumerate}
\vspace{10pt}
The numerator of $\Delta_e$ is divisible by the denominator due to the 
identity{\footnote{The identity is equivalent to \eqref{pppp}.}}
$$
\B_0(T) \B_1(-T) + \B_0(-T) \B_1(T) =2.
$$
Certainly, $\Delta_e$ is symmetric in the
half-edges.

\subsection{Pixton's relations ${\mathcal{P}}$}

%\begin{definition} \label{Not:relations}
Let $A = (a_1, \dots, a_n) \in \{0,1\}^n$.
Let
$$\cR_{g,A}^d\in \cS_{g,n}^d$$
 be the degree $d$ component of the strata algebra class 
$$
\sum_{\Gamma\in \mathsf{G}_{g,n}} \frac1{|{\text{Aut}}(\Gamma)| }
\, 
\frac1{2^{h^1(\Gamma)}}
\;
\left[\Gamma, \; \Bigl[
\prod \kappa_v \prod \B_l 
\prod \Delta_e
\Bigr ]_{\prod_v \zeta_v^{\mathrm{g}(v)-1}}
\right] \ \in \cS^*_{g,n},
$$
where the products are taken over all vertices, all legs, and all edges of the graph~$\Gamma$.
The subscript $\prod_v \zeta_v^{\mathrm{g}(v)-1}$ indicates
the coefficient of the monomial $\prod_v \zeta_v^{\mathrm{g}(v)-1}$
after the product inside the brackets is expanded. In fact,
$$\cR_{g,A}^d=0 \in \cS^d_{g,n}$$
unless the parity constraint
 $$g \equiv d+1+ \sum_{i=1}^n a_i \mod 2$$
holds.

We denote by $\tP$ the set of classes $\cR^d_{g,A}$ where
$$
d > \frac{g-1 + \sum_{i=1}^n a_i}{3}.
$$
By the following result, Pixton's set $\tP$ consists of tautological relations.

\begin{Theorem}[Janda 2013] \label{f5f}
Every element $\cR_{g,A}^d\in \tP$ lies in the kernel of 
the homomorphism $$q:\cS_{g,n}^* \rightarrow A^*(\oM_{g,n})\, .$$ 
%if
%$2g-2+n>0$, $a_i \in \{ 0, 1 \}$, and $d > \frac{g-1 + \sum_{i=1}^n a_i}{3}$.
\end{Theorem}

\vspace{2pt}
\subsection{Pixton's relations $\overline{\mathcal{P}}$} \label{q234}
%Pixton's set $\tP$ defines {\em finitely} many 
%tautological relations in each $\cS_{g,n}^d$.
The set $\tP$ is  extended to a larger set 
$$\tP \subset \overline{\mathcal{P}}$$
of tautological relations by the following construction.

The first step is to define a tautological
relation $\cR_{g,A,\sigma}^d\in \cS_{g,n}^d$ associated
to the data
\begin{enumerate}
\vspace{7pt}
\item[$\bullet$] $g,n\in \mathbb{Z}_{\geq 0}$ in the stable range $2g-2+n>0$,
\vspace{7pt}
\item[$\bullet$] $A=(a_1,\ldots, a_n), \ \,  a_i \in \mathbb{Z}_{\geq0}$,\ \, $a_i\equiv
 0 \text{\, or\, } 1 \mod 3$,
\vspace{7pt}
\item[$\bullet$] $\sigma$ is a partition of size $|\sigma|$ with 
parts  $\sigma_i \equiv  0$ or $1$ mod $3$,
\vspace{7pt}
\item[$\bullet$] $d\in \mathbb{Z}_{\geq 0}$ satisfying
$d > \frac{g-1+\sum_{i=1}^n a_i+|\sigma|}{3}$.
\end{enumerate}
\vspace{7pt}
Let $B=(b_1,\ldots,b_n,b_{n+1},\ldots,b_{n+\ell})$ be the unique vector
satisfying
$$b_j\in \{0,1\}\ \ \text{and}\ \
\begin{cases}
b_j\equiv a_j \ \ \ \mod 3 & 1\leq j \leq n \\
 b_{j} \equiv \sigma_{j-n} \mod 3& n+1\leq j \leq n+\ell 
\end{cases}\, ,$$
where $\ell$ is the length of $\sigma$.
Let $$\widehat{d}\, =\,  d- \sum_{j=1}^n  {\frac{a_j-b_j}{3}}
-\sum_{j=n+1}^{n+\ell} {\frac{\sigma_{j-n}-b_j}{3}} \, >\, 
\frac{g-1 + \sum_{j=1}^{n+\ell} b_j}{3}\, .$$
Pixton's definition{\footnote{Our conventions here
differ from \cite{Pix} by a global sign.}} of $\cR_{g,A,\sigma}^d\in \cS^d_{g,n}$ is 
$$ \cR_{g,A,\sigma} = p_{\ell*}\left(
\cR^{\widehat{d}}_{g,B} \, \cdot\, \prod_{j=1}^n \psi_j^{\frac{a_j-b_j}{3}}
\, \prod_{j=n+1}^{n+\ell} \psi_j^{1+\frac{\sigma_{j-n}-b_j}{3}} \right)\, ,$$
where $\cR^{\widehat{d}}_{g,B} \in \cS^{\widehat{d}}_{g,n+\ell}$
is in the set $\mathcal{P}$ and
$p_{\ell*}$ is push-forward by  the map forgetting the last $\ell$ markings,
$$p_{\ell*}: \cS^*_{g,n+\ell} \rightarrow \cS^*_{g,n}\, .$$
By Theorem \ref{f5f}, $\cR^{\widehat{d}}_{g,B}$ is a tautological relation.
Therefore,
$\cR_{g,A,\sigma}^d$ is also a
tautological relation. 

When $A=\emptyset$,  
the relations $\cR_{g,\sigma}^d\in \cS_{g,0}^d$ yield, after restriction
to $$\M_g \subset \oM_g\, ,$$
the Faber-Zagier relations of Theorem \ref{dddd}.{\footnote{The
relation $\cR_{g,\sigma,A}^d$ is trivial unless the parity condition
$$g \equiv d+|\sigma|+1+ \sum_{i=1}^n a_i \mod 2$$
holds.}} 

%$$ $$
%The $\mathbb{Q}$-linear span of $\tP$ is {\em not} an ideal in each
%$\cS_{g,n}$. Extend $\tP$ to the ideal generated by $I(\tP)$ in
%in each $\cS_{g,n}$. 
% For every forgetful map
%$$p_{m*}: \cS^{d+m}_{g,n+m} \rightarrow \cS^d_{g,n}$$
% but generates an ideal in each~$\cS_{g,n}$. 

Pixton's set $\overline{\mathcal{P}}$ is obtained
by taking the closure of the extended set of classes
\begin{equation}\label{xpxpxp} 
\left\{\, \cR^d_{g,A,\sigma}\in \cS_{g,n}^d\ \Big|\ d > \frac{g-1+\sum_{i=1}^n a_i+|\sigma|}{3}
\, \right\}
\end{equation}
under push-forward by all boundary maps:
add to the set \eqref{xpxpxp}  
 all classes
in $\cS_{g,n}^*$ which are obtained from a stable graph
$$\Gamma \in \mathsf{G}_{g,n}$$
with a class $\cR^{d(v)}_{g(v),A(v),\sigma(v)}$
placed on a single vertex $v\in \V(\Gamma)$ and any product of tautological 
classes placed on the other vertices of $\Gamma$.
By Theorem \ref{f5f}, every class in $\overline{\mathcal{P}}$
determines a tautological relation.

The subset ${\overline{\mathcal{P}}}$
lying in  a fixed $\cS^d_{g,n}$ is effectively computable. By the 
dimension restriction
$$d > \frac{g-1+\sum_{i=1}^n a_i+|\sigma|}{3}\, ,$$
only {\em finitely} many $\cR^d_{g,A,\sigma}$ lie in 
 ${\overline{\mathcal{P}}}$. The closure process by boundary
push-forward is
again {\em finite} because of the dimension restriction.
Hence, the $\mathbb{Q}$-linear span of ${\overline{\mathcal{P}}}$
in $\cS^d_{g,n}$ is generated by a finite list of
explicit classes.

Why stop at ${\overline{\mathcal{P}}}$? Why not consider
the closure with respect to further push-forwards and
pull-backs via the standard boundary and forgetful
maps? Pixton \cite{Pix,PixPhd}  has proven the set  ${\overline{\mathcal{P}}}$
is {\em closed} under all these further 
operations.{\footnote{From Pixton's results,
the most efficient definition of the $\mathbb{Q}$-linear
span of ${\overline{\mathcal{P}}}$
is as {\em the smallest set of ideals  
$$\{\, {\mathcal{I}}_{g,n} \subset \cS_{g,n}^*\, \}$$
which contains ${\mathcal{P}}$ and is closed under
the natural boundary and forgetful operations}.}}

\vspace{8pt}
\noindent {\bf Q6.} {\em Do  Pixton's relations ${\overline{\mathcal{P}}}$
span the ideal of relations among the tautological 
classes in $R^*(\oM_{g,n})$ for all $g$ and $n$?}
\vspace{8pt}

\subsection{Pixton's conjecture}
Pixton has conjectured an affirmative answer to question {\bf Q6}.
The evidence for Pixton's conjecture is (at least) the following:
\begin{enumerate}
\item[(i)] All the previously found relations 
occur in ${\overline{\mathcal{P}}}$.
The theory in genus $0$ is straightforward (and explained in \cite[Section 3.6]{PPZ}).
Modulo simpler relations, Getzler's genus 1 relation 
is
$$\cR^2_{1,(1,1,1,1)} \in \cS^2_{1,4}\, .$$
Pixton's conjecture is true{\footnote{The proof
uses the Gorenstein property of the
tautological rings $R^*(\oM_{0,n})$ and $R^*(\oM_{1,n})$.
The Gorenstein property is clear in $g=0$ since
$$R^*(\oM_{0,n})\cong H^*(\oM_{0,n})$$ and is proven in $g=1$ 
in \cite{PetG}.}} for $g\in\{0,1\}$ and
all $n$.
The genus 2 relation displayed in Section \ref{straa} is 
$$\cR^2_{2,(1,1,1)}\in \cS^2_{2,3}\, $$
modulo simpler relations.

\item[(ii)]  Computer calculations of Pixton's relations for low $(g,n)$ often
yield Gorenstein rings (forbidding further relations). However,
just as in the Faber-Zagier case, Pixton's relations do {\em not}
always yield Gorenstein rings.
\item[(iii)] Janda \cite{Janda15} has proven that a wide class
of semisimple Cohomological Field Theories (including
higher projective spaces and $r$-spin curves)
will not yield relations outside of $\overline{\mathcal{P}}$.
\end{enumerate}
Question {\bf Q6} has not been investigated
as extensively as question {\bf Q2} for the Faber-Zagier
relations since the moduli spaces of stable
curves are computationally more difficult to handle. 
However, the failure to find additional
$\kappa$ relations in $R^*(\M_g)$ may also be viewed
as supporting {\bf Q6} by the restriction property.

Pixton's proposal provides an effective calculus
of tautological classes on the moduli spaces $\oM_{g,n}$ of
stable curves. With an affirmative answer to {\bf Q6},
Pixton's calculus provides a {\em complete} answer.
Perhaps a reformulation of the set ${\overline{\mathcal{P}}}$
in a more directly algebraic setting will eventually be found.
How the subject will develop depends very much on the
answer to {\bf Q6}.  

Even a few years ago, a
calculus for the moduli space of curves seemed far out of reach.
Pixton's proposal has led to a striking change of outlook.

\subsection{Nonsingular and compact type curves} \label{ll99}
The moduli spaces
\begin{equation}\label{cfrt}
\M_{g,n}\, ,  \ \M_{g,n}^{\mathsf{rt}}\, ,  \ 
\M_{g,n}^{\mathsf{ct}} \  \subset\ \oM_{g,n}
\end{equation}
are all  open subsets. The algebras of tautological
classes{\footnote{See \cite{BSZ,Cav,Ion,JP,Tav1,Tav2} for various directions
in the study of these tautological rings.}}
 $$R^*(\M_{g,n})\subset A^*(\M_{g,n}), \  
R^*(\M_{g,n}^{\mathsf{rt}})\subset A^*( \M_{g,n}^{\mathsf{rt}}  ), \  
R^*(\M_{g,n}^{\mathsf{ct}})\subset A^*(\M_{g,n}^{\mathsf{ct}})
$$
are defined in each case as the image of the respective restriction of
$$R^*(\oM_{g,n}) \subset A^*(\oM_{g,n})\, .$$
A basic question here concerns the extension of
tautological relations over the boundary.

\vspace{8pt}
\noindent {\bf Q7.} {\em Does every tautological relation
in $R^*(\M_{g,n})$, $R^*(\M_{g,n}^{\mathsf{rt}})$, and
$R^*(\M_{g,n}^{\mathsf{ct}})$ arise from the restriction
of a tautological relation in $R^*(\oM_{g,n})$?
}
\vspace{4pt}

If the answers to  {\bf Q6} {\em and} {\bf Q7} are both affirmative, then
the Pixton calculus determines the tautological
rings in all the nonsingular and compact type cases \eqref{cfrt}.

\vspace{6pt}
The $n^{th}$ fiber product $\mathcal{C}^n_g$ 
of the universal
curve 
$$\pi:\mathcal{C}_g \rightarrow \M_g$$ is
 {\em not} an open set of $\oM_{g,n}$. However,
there is a proper surjection
$$\M_{g,n}^{\mathsf{rt}} \stackrel{\epsilon}{\longrightarrow}  \mathcal{C}^n_g\, .$$
The tautological ring $R^*(\mathcal{C}^n_g)$ is defined{\footnote{For
another definition of $R^*(\mathcal{C}^n_g)$ and further study of
the relationship with
$R^*(\M_{g,n}^{\mathsf{rt}})$, see \cite{PetP}.}}
as the image of $R^*(\M_{g,n}^{\mathsf{rt}})$ under $\epsilon_*$.
If the answers to {\bf Q6} and {\bf Q7} are both affirmative,
then Pixton's calculus also determines $R^*(\mathcal{C}^n_g)$.

\subsection{Further directions}

\subsubsection{Symmetries}
The symmetric group $\Sigma_n$ acts naturally on $\oM_{g,n}$
by permuting the $n$ markings. Since both $\cS_{g,n}^*$ and $R^*(\oM_{g,n})$
carry induced  $\Sigma_n$-representations and  
$$q: \cS_{g,n}^* \rightarrow R^*(\oM_{g,n})$$
is a morphism of 
$\Sigma_n$-representations, the ideal of
tautological relations also carries an induced $\Sigma_n$-representation.

Pixton's relations interact in interesting
ways with the $\Sigma_n$-action. For $g>0$,
Pixton proves \cite[Proposition 2]{Pix}
that
the {\em new} relations{\footnote{The {\em new} relations
are those which do not lie in the $\mathbb{Q}$-span of
relations coming from lower genus,  lower marking number,
or lower codimension.}}
in $\overline{\mathcal{P}}^d_{g,n}$ are generated by 
$\Sigma_n$-invariant tautological relations. 
For example,
both Getzler's relation
and the genus 2 relation displayed in Section \ref{straa} are new
and invariant.

The symmetric group acts on the entire cohomology $H^*(\oM_{g,n})$.
For $g\leq 2$, the symmetric group representations are well 
understood \cite{vdgF2,vdgF1,vdgF11, Get0, Get1,
PetAb}, and there is significant
progress \cite{vdgF22} in genus $3$. Constraints
on the $\Sigma_n$-action on $R^*(\oM_{g,n})$ are proven in
\cite[Section 4]{FP3} and show certain cohomology classes can {\em not} be
tautological.

\subsubsection{Push-forward relations}
A tautological relation in $\cS^*_{g,n}$ yields a universal
equation for the genus $\leq g$ Gromov-Witten theory
 of {\em every}
target variety $X$ \cite{KKMM}. However, in genus 1, such tautological
relations do not appear to be enough to prove the Virasoro
constraints  \cite{EHX} for arbitrary targets \cite{Xliu}. Can any further
universal equations in Gromov-Witten theory be 
found in the geometry of the
moduli spaces of curves?

An idea to find further universal equations in
 Gromov-Witten theory using push-forwards is the following.
Consider the gluing map, 
$$\delta: \oM_{g,2} \rightarrow \oM_{g+1}\, ,$$
with image equal to the divisor $\Delta_0\subset \oM_g$ of
curves with a nonseparating node. Elements of 
the kernel
$$\delta_*: R^*(\oM_{g,2}) \rightarrow R^{*+1}(\oM_{g+1})$$
yield  universal equations for the genus $\leq g+1$ Gromov-Witten theory
of targets $X$ with {\em no genus $g+1$ terms}. Hence, 
elements of $\text{ker}(\delta_*)$ yields universal
equations for genus $\leq g$ Gromov-Witten theory.

To formalize the notion of a {push-forward relation}, consider
the composition
 $$\cS_{g,2}^* \stackrel{q}{\longrightarrow}
R^*(\oM_{g,2}) \stackrel{\delta_*}{\longrightarrow}  R^{*+1}(\oM_{g+1})\, .$$
We have $\text{ker}(q) \subset \text{ker}(\delta_*\circ q)$.
The gluing map $\delta$ lifts to
$$\widetilde{\delta}: \cS_{g,2}^* \rightarrow \cS^{*+1}_{g+1}\, ,$$
and we also have $\text{ker}(\widetilde{\delta}) \subset \text{ker}(\delta_*\circ q)$.
The subspace 
$${\ker(q)+
\ker(\widetilde{\delta})} \subset {\text{ker}(\delta_*\circ q)}$$
yields relations in Gromov-Witten theory already
captured before considering $\delta$.

A {\em push-forward relation} is a nonzero element
of $$\frac{\text{ker}(\delta_*\circ q)}{\ker(q)+
\ker(\widetilde{\delta})}\, .$$
Push-forward relations yield universal contraints
in Gromov-Witten theory which appear to go beyond
the contraints obtained from tautological relations in
$\cS^*_{g,n}$.

Are there any push-forward relations? Possible
candidates were found in \cite{LiuP}.
Let $r$ be even and satisfy $2\leq r\leq g-1$. Define
 $$\chi_{g,r}\, =\, \sum_{a+b=2g+r} (-1)^a [\Gamma_\bullet,\psi_1^a \psi_2^b]\ \in \cS_{g,2}^{2g+r}. $$
Here, $\Gamma_\bullet$ is the unique graph with no edges.
By \cite[Theorem 2]{LiuP}, 
$$\delta_* q (\chi_{g,r})=0\, .$$
Moreover,
%A straightforward argument{\footnote{Suppose there exists an expression
%$$\chi_{g,r}= \alpha + \beta \, , \ \ \ \alpha \in \text{ker}(q)\,, \ \ \beta\in
%\text{ker}(\widetilde{\xi})\, .$$
%Let $\sigma\in \Sigma_2$ be the 2-cycle. By the 
%$\Sigma_2$-invariance of $\chi_{g,r}$,
%$$2 \chi_{g,r}= \alpha+ \sigma(\alpha)+ \beta+\sigma(\beta).$$
%Since every element $\beta\in \text{ker}(\widetilde{\xi})$ satisfies
%$$\beta+\sigma(\beta)=0 \in \text{ker}(\widetilde{\xi})\, ,$$
%we can then conclude $q(\chi_{g,r})= q(\alpha+\sigma(\alpha))=0$ 
%which contradicts the nonvanishing integrals of 
%\cite[Section 0.3]{LiuP}.}}
the results of \cite{LiuP}
show $$\chi_{g,r} \ \notin\  \text{ker}(q) \ \ \text{and} \ \ 
\chi_{g,r} \ \notin\ \text{ker}(\widetilde{\delta})\, ,$$
but whether $\chi_{g,r}$ ever avoids the
sum  $\text{ker}(q) +\text{ker}(\widetilde{\delta})$ is open.

Pixton's relations can be used to search for push-forward
relations (via several related constructions) which have the possibility of
producing new universal Gromov-Witten
equations. What role such relations will play in Gromov-Witten
theory is not yet known.

\section{Double ramification cycles}

\subsection{Overview}
Curves of genus $g$ which admit a map to $\mathbb{P}^1$
with specified ramification profile $\mu$ over $0\in \mathbb{P}^1$
and $\nu$ over $\infty \in \mathbb{P}^1$
define a double ramification cycle $\DR_{g}(\mu,\nu)$ on the
moduli space of curves.
The restriction of the double ramification cycle to the moduli space of
nonsingular curves is a classical
topic related to the linear equivalence of divisors.

The cycle $\DR_{g}(\mu,\nu)$ on the moduli space of stable curves
is defined via the virtual fundamental class of the moduli space of
stable
maps to rubber. 
An explicit formula for $\DR_{g}(\mu,\nu)$
in the tautological ring, conjectured by Pixton \cite{PixDR} in 2014
and proven in \cite{jppz} in 2015, 
is presented here. 
Pixton's double ramification formula
expresses the cycle as a sum over stable graphs
(corresponding to strata classes) with summands given by a
product of leg and edge factors.
The result shows how the calculus of
tautological classes works in practice.

\subsection{Moduli of relative stable maps}
Let ${\mu}=({\mu_1}, \ldots, {\mu_{\ell(\mu)}})$ and
${\nu}=({\nu_1}, \ldots,
{\nu_{\ell(\nu)}})$
be partitions of equal size,
$$\sum_{{i=1}}^{\ell(\mu)} {\mu_i} = \sum_{{j=1}}^{\ell(\nu)} 
{\nu_j}\, .$$
Let $C$ be a genus $g$ curve.
We consider maps 
$$f:C \rightarrow \mathbb{P}^1$$
with
ramification profiles ${\mu}$ over $0\in \mathbb{P}^1$ and
and ${\nu}$ over $\infty\in {\mathbb{P}^1}$.
Two such maps 
$$C\stackrel{f}{\longrightarrow}\mathbb{P}^1\, , \ \ \
C\stackrel{\widehat{f}}{\longrightarrow}\mathbb{P}^1$$
are declared equivalent if $f$ and $\widehat{f}$
differ by a reparameterization
of the target which keeps both $0$ and $\infty \in \mathbb{P}^1$ fixed.

{\centering{\ \ \ \ \hspace{50pt} \includegraphics[width=6.3cm]{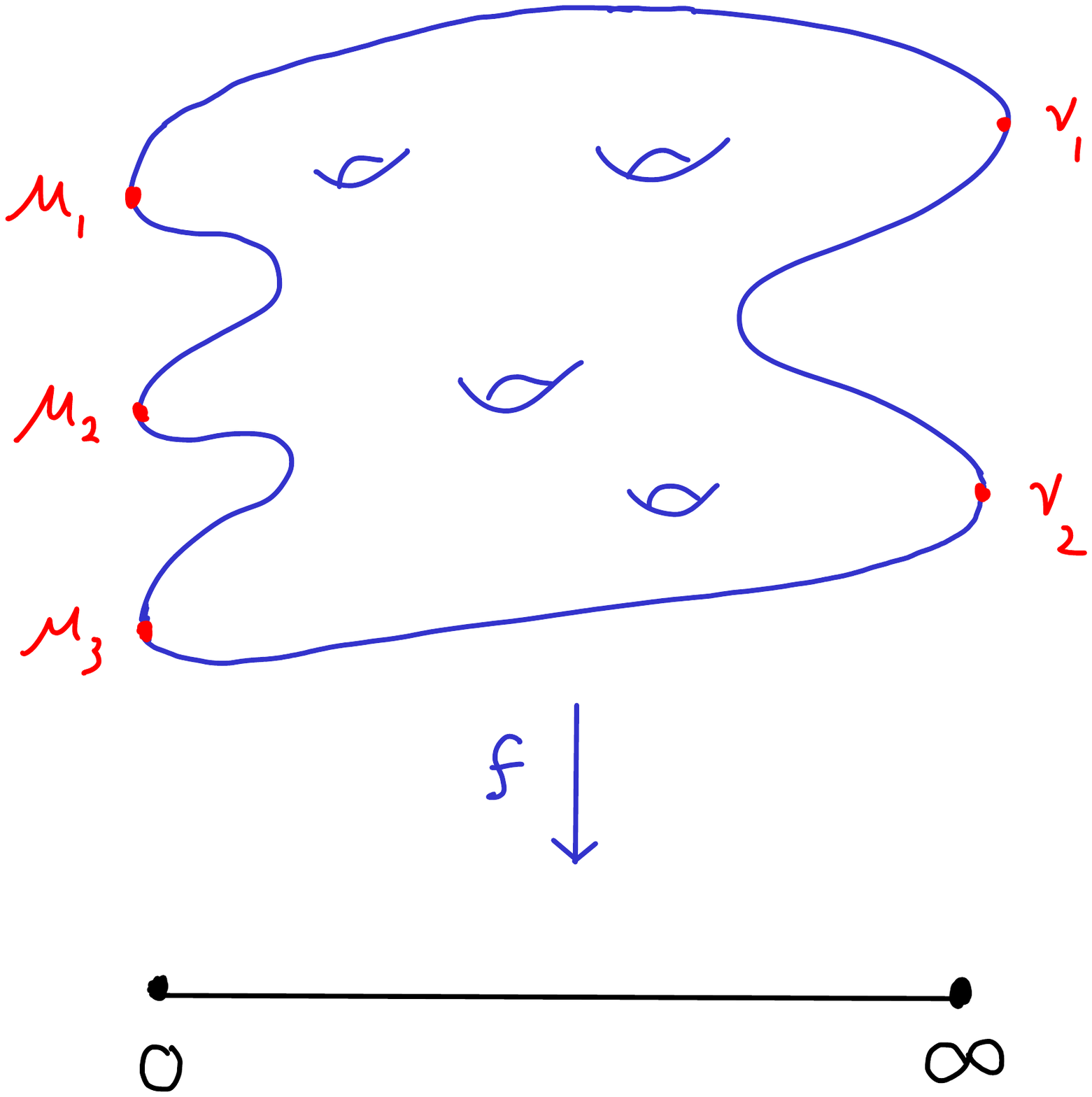}}}

\vspace{-40pt}

A natural compact moduli space of such maps $f$ arises in 
Gromov-Witten theory. Let
$$\oM_g({\mathbb{P}}^1,\mu,\nu)^\sim$$
be the moduli space of {\em stable relative maps} to rubber
with ramification profiles $\mu$ and $\nu$. 
In the moduli of relative stable maps,
 $f$ may degenerate in several ways: the domain $C$ may acquire
nodes, $f$ may be constant on irreducible components of $C$, and
the target $\mathbb{P}^1$ may degenerate. The first two phenomena are
illustrated in the following diagram.

{\centering{\ \ \ \ \hspace{50pt} \includegraphics[width=6.3cm]{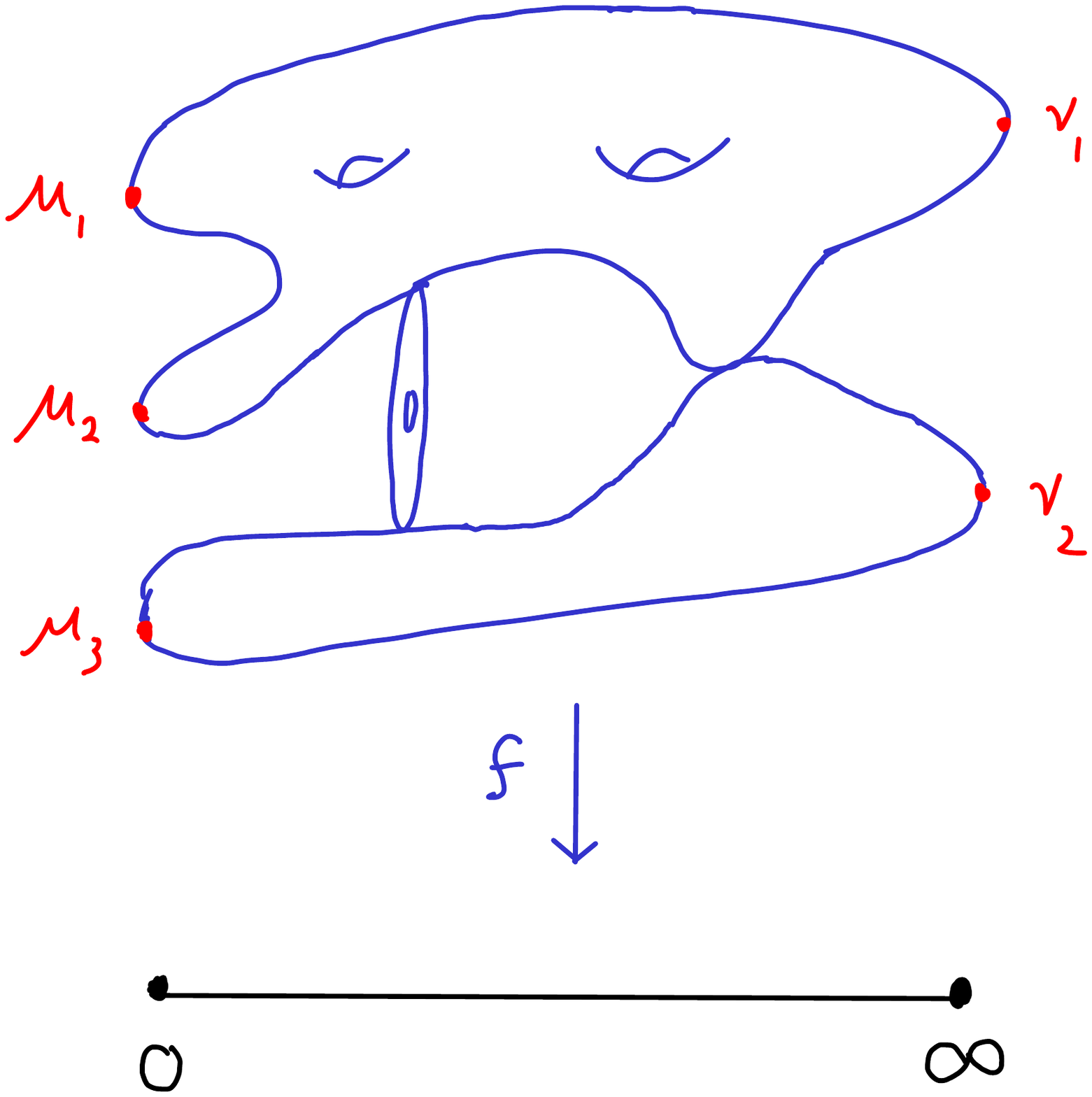}}}

\vspace{-30pt}
\noindent 
We refer the reader 
to \cite{LR,JL1,JL2} for the basic definitions 
of relative Gromov-Witten theory  and the  foundational development.

\subsection{Double ramification cycles}
There is a natural morphism
$$ \rho: {\oM_{g}}({\mathbb{P}^1},{\mu},
{\nu})^{\sim} 
 \rightarrow {\oM_{g,{\ell(\mu)}
+ {\ell(\nu)}}}$$
forgetting everything except the {marked domain curve}.
The {\em double ramification cycle} is the push-forward 
of the {virtual fundamental class}{\footnote{The 
expected dimension of  ${\oM_{g}}({\mathbb{P}^1},
{\mu},
{\nu})^{\sim}$ is $3g-3+\ell(\mu)+\ell(\nu)-g$ where
$g$ is the dimension of the Jacobian of the domain curve.}}
$${\mathsf{DR}_{{g}}({\mu},{\nu}})\  =\ 
\rho_*\Big[{\oM_{g}}({\mathbb{P}^1},
{\mu},
{\nu})^{\sim}\Big]^{{vir}} \ \in 
{A}^{g}({\oM_{g,{\ell(\mu)}
+ {\ell(\nu)}        }})\ .$$

\noindent
Eliashberg asked in 2001: {\em what is the formula for
${\mathsf{DR}_{{g}}({\mu},{\nu}})$}?
As a first step,
the double ramification cycle was proven to lie in
the tautological ring in \cite{FPrel} in 2005,
$${\mathsf{DR}_{{g}}({\mu},{\nu}})\in R^g(\oM_{g,\ell(\mu)+\ell(\nu)})\, .$$
The restriction of ${\mathsf{DR}_{{g}}({\mu},{\nu}})$
to the moduli space $\M_{g,\ell(\mu)+\ell(\nu)}^{\mathsf{ct}}$
of curves of compact type can be calculated via the
geometry of the universal Jacobian.{\footnote{The matching
on the moduli space of curves of compact type
of the definitions of the double ramification cycle 
via the virtual class and the Jacobian geometry is
not trivial and is proven in \cite{Cav2,MW}.}}
 The result is
Hain's formula \cite{GruZak, Hain}.

\pagebreak
\subsection{Pixton's formula}
\subsubsection{Ramification vector}
We place the ramification data in a vector
$$
({\mu_1},{\ldots},{\mu_{\ell(\mu)}}, {-\nu_1}, {\ldots}, {-\nu_{\ell(\nu)}})\, .
$$
For any vector ${S}=({s_1},\ldots,{s_n})$ with $\sum_{{i}} {s_i}=0$, we have 
$${\mathsf{DR}_{{g},{S}}}\  \in 
{R}^{g}({\oM_{g,{n}}})\ .$$
The positive parts of $S$ specify ramification over $0\in \mathbb{P}^1$ and
the negative parts specify ramification over $\infty \in \mathbb{P}^1$.
Free points corresponding to part of $S$ equal to $0$ are also permitted.

{Pixton} \cite{PixDR} conjectured a beautiful formula 
for ${\mathsf{DR}_{{g},{S}}}\in R^g(\oM_{g,n})$ which involves
a sum over admissible weightings of stable graphs.

\subsubsection{Admissible weightings}
Let ${S}=({s_1},{\ldots},{s_n})$ be  double ramification data. 
Let ${\Gamma}\in \mathsf{G}_{g,n}$ be a {stable graph} of genus ${g}$ with ${n}$ legs.
An {\em {admissible weighting}} is 
a function on the set of half-edges,
$$ w:{\H}({\Gamma}) \rightarrow \mathbb{Z},$$
which satisfies: 
\vspace{5pt}
\begin{enumerate}
\item[(i)] $\forall {h_i}\in {\L}({\Gamma})$, 
$w({h_i})={s_i}$,
\item[(ii)] $\forall {e} \in {\E}({\Gamma})$ consisting
of the half-edges 
${h}({e}),{h'}({e}) \in {\H}({\Gamma})$,
$$w({h})+w({h'})=0\, ,$$
\item[(iii)] $\forall {v}\in {\V}({\Gamma})$,
$\sum_{{v}({h})= {v}} w({h})=0$. 
\end{enumerate}
\vspace{5pt}
A {stable graph $\Gamma$}, however, may have {\em infinitely} many 
{admissible weightings} $w$. 

In order to regularize the
sum over admissible weightings, Pixton introduced a regularization parameter
$${r}\in \mathbb{Z}_{>0}\, .$$
An {\em {admissible weighting} mod ${r}$}
of ${\Gamma}$ is a function,
$$ w:{\H}({\Gamma}) \rightarrow \{{0},{\ldots},{r-1}\},$$
which satisfies{\footnote{For example, for (i), we require
$w({h_i})\equiv{s_i} \mod {r}$.}} 
the conditions of (i-iii) above $mod$ ${r}$.
Let $\mathsf{W}_{{\Gamma},{r}}$ be the set of admissible weightings mod 
${r}$
of ${\Gamma}$. 
The set $\mathsf{W}_{{\Gamma},{r}}$ is {\em finite} and can be summed
over.

\subsubsection{Formula}

Let ${r}$ be a positive integer.
We denote by
$${\mathsf{Q}_{{g},{S}}^{d,{r}}}\in {R}^d({\oM_{g,n}})$$ 
the degree ${d}$ component of 
the class 
\begin{multline*}
\hspace{-10pt}\sum_{{\Gamma}\in {\mathsf{G}_{g,{n}}}} 
\sum_{w\in \mathsf{W}_{{\Gamma},{r}}}
\frac1{|\text{Aut}({\Gamma})| }
\, 
\frac1{r^{h^1({\Gamma})}}
\;
\xi_{{\Gamma}*}\Bigg[
\prod_{i=1}^n \exp({s_i}^2 {\psi}_{{h_i}}) \cdot 
\\ \hspace{+10pt}
\prod_{{e}=({h},{h'})\in {\V}({\Gamma})}
\frac{1-\exp(-w({h})w({h'})({\psi}_{{h}}+{\psi}_{{h'}}))}
{{\psi}_{{h}} + {\psi}_{{h'}}} \Bigg]\, .
\end{multline*}
For fixed ${g}$, ${S}$, and ${d}$, the 
class
${\mathsf{Q}_{{g},S}^{d,r}} \in {R}^{{d}}({\oM_{g,n}})$
is polynomial in ${r}$ for sufficiently large ${r}$, see
\cite[Appendix]{jppz}.
We denote by ${\mathsf{P}_{{g},S}^d}$ the value at ${r}={0}$, 
$${\mathsf{P}_{{g},S}^d}\, =\, {\mathsf{Q}_{{g},S}^{d,r}}\big|_{r=0}\, .$$ 
Hence, ${\mathsf{P}_{{g},S}^d}$ is the {\em {constant}} term.

%{Pixton} conjectured in \bbb{2014} the following result.

\begin{Theorem} [Janda-P.-Pixton-Zvonkine 2015]\label{jppz3} For all $g\geq 0$
and double ramification data $S$,
$${\mathsf{DR}_{{g},S}} = 2^{-{g}}\, 
{\mathsf{P}_{{g},S}^{{g}}}\, \in {R}^{{g}}
({\oM_{g,n}})\, .$$ 
\end{Theorem}

\vspace{2pt}
\noindent 
When restricted to $\M_{g,n}^{\mathsf{ct}}$, Theorem \ref{jppz3}
 recovers Hain's formula.

The proof of Theorem \ref{jppz3} uses the equivariant {Gromov-Witten} theory of 
 ${\mathbb{P}}^1$ with:
\begin{enumerate}
\item[$\bullet$] an orbifold ${B\mathbb{Z}_r}$-point at 
$0\in {\mathbb{P}}^1$,
\item[$\bullet$] a relative point at
 $\infty \in {\mathbb{P}}^1$.
\end{enumerate}
Hence, {orbifold GW} theory \cite{AGV,ChenR}, {relative GW} theory \cite{LR,JL1,JL2}, 
and the virtual localization formula \cite{GP} all play a role.
Over ${B\mathbb{Z}_r}$, Hurwitz-Hodge integrals arise exactly
in the form of \cite{JPT} and are analysed via Chiodo's formula
\cite{Chiodo}. The double ramification cycle arises over the relative point
$\infty\in \mathbb{P}^1$.

In addition to conjecturing the formula of Theorem \ref{jppz3}, Pixton \cite{PixDR}
conjectured the following vanishing proven in \cite{cj}.

\begin{Theorem}[Clader-Janda 2015] \label{kqkq}
For all $g\geq 0$, double ramification data $S$, and $d>g$, 
$${\mathsf{P}}_{g,S}^d = 0 \, \in R^d(\oM_{g,n}) .$$ 
\end{Theorem}

\noindent Clader and Janda further prove that the tautological relations obtained
from Theorem \ref{kqkq} are a consequence of Pixton's set $\overline{\mathcal{P}}$
discussed in Section \ref{q234}.
For $d<g$,  the classes ${\mathsf{P}}_{g,S}^d$ do not yet have a
geometric interpretation. 

\subsection{The Hodge bundle}
Let $g\ge2 $ and $S=\emptyset$, so $\mu=\nu=\emptyset$.
Then, the morphism $\rho$ from  
the moduli of stable maps to the moduli of curves is 
an {\em isomorphism},
$$ \rho: \oM_{g}(\mathbb{P}^1,\emptyset,\emptyset)^{\sim} 
\, \stackrel{\sim}{\longrightarrow} \,
 \oM_{g}\ .$$
By a study of the obstruction theory,
$$\DR_{g,\emptyset} = (-1)^g\lambda_g \ \in R^g(\oM_{g})\, , $$
 where $\lambda_g$ is the top Chern class of the Hodge bundle
$$\mathbb{E} \rightarrow \oM_{g}\ . $$

Pixton's formula in the $S=\emptyset$ case therefore yields an expression
for $\lambda_g$.
By the analysis of \cite[Section 0.5.3]{jppz}, the result is
a new and very special formula for $\lambda_g$: all
the strata which appear are supported on 
$$\Delta_0 \subset \oM_g\, ,$$
the divisor with a nonseparating node.{\footnote{In particular, the
vanishing
$$\lambda_g^2=0 \in R^*(\oM_g)$$
is an immediate consequence since $\lambda_g|_{\Delta_0}=0$.}}

In the diagrams{\footnote{The diagrams (taken from \cite{jppz}) are 
computed using code written by A. Pixton. The
artistic display is due to F. Janda.}} below,
each labeled graph $\Gamma$ describes a moduli space $\oM_\Gamma$, 
a tautological class $\gamma$, and 
a natural map 
$$\xi_\Gamma: \oM_\Gamma \to 
\oM_g\, .$$ 
The convention in the diagrams is that the labeled graph 
represents the 
cycle class 
$$\xi_{\Gamma*}(\gamma)=q\left([\Gamma,\gamma]\right)\, ,$$
 see Section \ref{straa}.

\tikz{\coordinate (A) at (0,0); \coordinate (B) at (1,0); \coordinate (C) at (0.6,0.5);}
\tikzset{baseline=0, label distance=-3mm}
\def\NC{\draw (0,0.25) circle(0.25);}
\def\NL{\draw plot [smooth,tension=1.5] coordinates {(0,0) (-0.2,0.5) (-0.5,0.2) (0,0)};}
\def\NR{\draw plot [smooth,tension=1.5] coordinates {(0,0) (0.2,0.5) (0.5,0.2) (0,0)};}
\def\NN{\NL\NR}
\def\NNN{\NN \begin{scope}[rotate=180] \NR \end{scope}}
\def\NNNN{\NN \begin{scope}[rotate=180] \NN \end{scope}}
\def\NRS{\begin{scope}[shift={(B)}] \NR \end{scope}}
\def\NRD{\begin{scope}[rotate around={-90:(B)}] \NRS \end{scope}}
\def\DE{\draw plot [smooth,tension=1] coordinates {(0,0) (0.5,0.15) (1,0)}; \draw plot [smooth,tension=1.5] coordinates {(0,0) (0.5,-0.15) (1,0)};}
\def\DES{\begin{scope}[shift={(B)}] \DE \end{scope}}
\def\TE{\DE \draw (A) -- (B);}
\def\QE{\DE \draw plot [smooth,tension=1] coordinates {(A) (0.5,0.05) (B)}; \draw plot [smooth,tension=1.5] coordinates {(A) (0.5,-0.05) (B)};}
\def\T{\draw (0.2,0) -- (C) -- (B) -- (0.2,0);}
\def\TT{\draw (C) -- (B) -- (0.2,0); \draw plot [smooth,tension=1] coordinates {(0.2,0) (0.3,0.3) (C)}; \draw plot [smooth,tension=1] coordinates {(0.2,0) (0.5,0.2) (C)};}
\newcommand{\nn}[3]{\draw (#1)++(#2:3mm) node[fill=white,fill opacity=.85,inner sep=0mm,text=black,text opacity=1] {$\substack{\psi^#3}$};}
\renewcommand{\gg}[2]{\fill (#2) circle(1.3mm) node {\color{white}$\substack #1$};}

%\paragraph{Genus~1.}
%\begin{equation*}
% \lambda_1 = \frac 1{24} \tikz{\NC \gg{0}{A}}.
%\end{equation*}

\pagebreak
\noindent The formulas for $\lambda_2\in R^2(\oM_2)$
and $\lambda_3\in R^3(\oM_3)$ are: 

%\vspace{-5pt}
\begin{equation*}
\lambda_2 = 
  \frac 1{240} \tikz{\NC \nn{A}{130}{{}} \gg{1}{A}}
  + \frac 1{1152} \tikz{\NN \gg{0}{A}},
\end{equation*}

%\vspace{-20pt}
\begin{align*}
\lambda_3 &= 
  \frac 1{2016} \tikz{\NC \nn{A}{130}{2} \gg{2}{A}}
  + \frac 1{2016} \tikz{\NC \nn{A}{130}{{}} \nn{A}{50}{{}} \gg{3}{A}}
  - \frac 1{672} \tikz{\DE \nn{A}{30}{{}} \gg{1}{A} \gg{1}{B}}
  + \frac 1{5760} \tikz{\NN \nn{A}{160}{{}} \gg{1}{A}} \\
  &
  - \frac{13}{30240} \tikz{\TE \gg{0}{A} \gg{1}{B}}
  - \frac 1{5760} \tikz{\NL \DE \gg{0}{A} \gg{1}{B}}
  + \frac 1{82944} \tikz{\NNN \gg{0}{A}}.
\end{align*}

\noindent More interesting is the formula for $\lambda_4\in R^4(\oM_4)$:

$$
  \frac 1{11520} \tikz{\NC \nn{A}{130}{3} \gg{3}{A}}
  + \frac 1{3840} \tikz{\NC \nn{A}{130}{2} \nn{A}{50}{{}} \gg{3}{A}}
  - \frac 1{2880} \tikz{\DE \nn{A}{30}{2} \gg{1}{A} \gg{2}{B}}
  - \frac 1{3840} \tikz{\DE \nn{A}{30}{{}} \nn{A}{-30}{{}} \gg{1}{A} \gg{2}{B}}
  - \frac 1{1440} \tikz{\DE \nn{A}{30}{{}} \nn{B}{150}{{}} \gg{1}{A} \gg{2}{B}} 
$$
$$
  - \frac 1{1920} \tikz{\DE \nn{A}{30}{{}} \nn{B}{-150}{{}} \gg{1}{A} \gg{2}{B}}
  - \frac 1{2880} \tikz{\DE \nn{B}{150}{2} \gg{1}{A} \gg{2}{B}}
  - \frac 1{3840} \tikz{\DE \nn{B}{150}{{}} \nn{B}{-150}{{}} \gg{1}{A} \gg{2}{B}}
  + \frac 1{48384} \tikz{\NN \nn{A}{160}{2} \gg{2}{A}}
$$
$$
  + \frac 1{48384} \tikz{\NN \nn{A}{160}{{}} \nn{A}{110}{{}} \gg{2}{A}} 
  + \frac 1{115200} \tikz{\NN \nn{A}{160}{{}} \nn{A}{20}{{}} \gg{2}{A}}
  + \frac 1{960} \tikz{\T \nn{B}{180}{{}} \gg{1}{B} \gg{1}{0.2,0} \gg{1}{C}}
  - \frac{23}{100800} \tikz{\TE \nn{A}{-30}{{}} \gg{2}{A} \gg{0}{B}}
$$
$$
  - \frac 1{57600} \tikz{\DE \NRS \nn{B}{20}{{}} \gg{2}{A} \gg{0}{B}} 
  - \frac 1{16128} \tikz{\DE \NRS \nn{B}{-150}{{}} \gg{2}{A} \gg{0}{B}}
  - \frac 1{16128} \tikz{\DE \NRS \nn{A}{-30}{{}} \gg{2}{A} \gg{0}{B}}
$$
$$
  - \frac 1{57600} \tikz{\DE \NRS \nn{B}{20}{{}} \gg{1}{A} \gg{1}{B}}
  - \frac 1{16128} \tikz{\DE \NRS \nn{B}{-150}{{}} \gg{1}{A} \gg{1}{B}} 
   - \frac 1{16128} \tikz{\DE \NRS \nn{A}{-30}{{}} \gg{1}{A} \gg{1}{B}}
  - \frac{23}{100800} \tikz{\TE \nn{A}{-30}{{}} \gg{1}{A} \gg{1}{B}}
$$
$$
  + \frac{23}{100800} \tikz{\TT \gg{2}{B} \gg{0}{0.2,0} \gg{0}{C}}
  + \frac{23}{50400} \tikz{\TT \gg{1}{B} \gg{1}{0.2,0} \gg{0}{C}}
  + \frac 1{16128} \tikz{\T \NRS \gg{0}{B} \gg{1}{0.2,0} \gg{1}{C}} 
  + \frac 1{115200} \tikz{\DE \DES \gg{1}{A} \gg{0}{B} \gg{1}{2,0}}
$$
$$
  + \frac 1{276480} \tikz{\NNN \nn{A}{20}{{}} \gg{1}{A}}
  - \frac{13}{725760} \tikz{\NL \TE \gg{1}{A} \gg{0}{B}}
  - \frac 1{138240} \tikz{\NL \NRS \DE \gg{1}{A} \gg{0}{B}} 
$$
$$
  - \frac{43}{1612800} \tikz{\QE \gg{1}{A} \gg{0}{B}}
  - \frac{13}{725760} \tikz{\NRS \TE \gg{1}{A} \gg{0}{B}}
$$
$$
  - \frac 1{276480} \tikz{\NRS \NRD \DE \gg{1}{A} \gg{0}{B}}
  - \frac 1{7962624} \tikz{\NNNN \gg{0}{A}}\, .
$$

\vspace{5pt}

\subsection{Further formulas}
Pixton's formula for the double ramification cycle was presented
here as an example. Several other formulas have been recently
studied (the Chern character of the Verlinde bundle \cite{MOPPZ},
the cycle class of the loci of holomorphic/meromorphic differentials \cite[Appendix]{FarP}).
The form of a summation over stable graphs $\mathsf{G}_{g,n}$ 
with summands given by a product over vertex, leg, and edge
factors is ubiquitous (and reminiscent of Feynman expansions of
integrals in quantum field theory).

\vspace{+20 pt}
\noindent Departement Mathematik \\
\noindent ETH Z\"urich \\
\noindent R\"amistrasse 101\\
\noindent 8092 Z\"urich,  Switzerland

\end{document}